\documentclass[11pt]{article}

\newtheorem{theorem}{Theorem}
\newtheorem{lemma}{Lemma}
\newtheorem{proof}{Proof}
\usepackage{hyperref}

\usepackage{amsmath}
\usepackage{graphicx}
\usepackage{amssymb}
\usepackage{graphicx}
\usepackage[dvipsnames,usenames]{color}
\usepackage{subfigure}
\input{epsf.tex}

\usepackage[affil-it]{authblk}
\title{A two-class queueing system with constant retrial policy and general class dependent service times}
\author{Ioannis Dimitriou\footnote{ idimit@math.upatras.gr}}
\affil{\small Department of Mathematics, 
University of Patras, P.O.~Box 
26500, Patras, Greece.}
\date{\small \today}
\begin{document}
\maketitle
\begin{abstract}
A single server retrial queueing system with two-classes of orbiting customers, and general class dependent service times is considered. If an arriving customer finds the server unavailable, it enters a virtual queue, called the orbit, according to its type. The customers from the orbits retry independently to access the server according to the constant retrial policy. We derive the generating function of the stationary distribution of the number of orbiting customers at service completion epochs in terms of the solution of a Riemann boundary value problem. For the symmetrical system we also derived explicit expressions for the expected delay in an orbit without solving a boundary value problem. A simple numerical example is obtained to illustrate the system's performance.
\end{abstract}

\section{Introduction}
Queueing systems with retrial customers are characterized by the feature that an arriving customer who finds the server unavailable, departs temporarily from the system, and repeats its attempt to connect with the server after some random time according to a specific access policy. The so called repeated customers are temporarily stored in a pool of unsatisfied customers (called orbit or retrial group), and are superimposed on the normal stream of external arrivals. For a complete review of the main results, the interested reader is referred to the seminal books \cite{fal,art}, and in the detailed review papers \cite{art2,tuan1,kim}.
\subsection{Related work and applications}
Single class retrial systems under constant retrial policy were investigated in \cite{artw,avr2,choi1,choi2,du,far,gao,wang} (not exhaustive list). Clearly, there have been very limited results in retrial queueing literature with multiple classes of retrial customers. A two class retrial system with arbitrary distributed service requirements and classical retrial policy was firstly analyzed in \cite{kulk}, whereas the extension to an arbitrary number of classes of retrial customers was investigated in \cite{fali}. In \cite{mou} a non-preemptive priority mechanism was included in the work in \cite{kulk,fali}, while in \cite{land} a multiclass retrial queue with many phases of service was also investigated. In all the above mentioned works, a classical retrial policy was used and the authors derived expressions for the expected number of customers in orbit queues. Recently, there has been a lot of attention to the application of polling retrial systems with glue periods on the modeling of optical networks \cite{abi1,abi2,abi3,box}. In \cite{avr}, the authors studied a two-class system with common exponential service requirements and constant retrial policy. Their analysis led to a functional equation, which is solved with the aid of the theory of Riemann-Hilbert boundary value problems. Several generalizations of this model by considering coupled orbit queues, and simultaneous arrivals were considered in \cite{dim1,dim2}. A two class retrial system with common arbitrarily distributed paired service, and potential applications in wireless systems under network coding was investigated in \cite{dim}.

In general, multiclass retrial systems with constant retrial policy serve as a model for competing job streams in a carrier sensing multiple access system, where the jobs, after a failed attempt to network access, wait in an orbit queue; e.g., a local area computer network with bus architecture where the different types of customers can be interpreted as customers with different priority requirements \cite{spa}. Under the constant retrial policy we are able to stabilize and control the multiple access system. Such a priority setting can also be applied to train or vehicular onboard networks. In such a case the high priority jobs correspond to critical system control signals, and the low priority jobs correspond to onboard passenger internet access traffic. 

Other potential applications may be found in cooperative wireless systems. Such systems consist of a finite number of source users that transmit packets to a common destination node, and a finite number of assistant users, called relay nodes (i.e., the orbit queues) that assist them by retransmitting their failed packets; e.g., \cite{pap,pap2,dim1,dim2}. Other applications can be found in telecommunication systems with call-back option in call centers \cite{du,tu}, where an operator (i.e. a server) calls-back an unsatisfied customer after some random time.
\subsection{Our contribution}
The important feature of this work is the two class setting under constant retrial policy, and arbitrarily distributed service requirements, which depend on the type of the job as well as the instant of its arrival. In particular, the service times of primary jobs that occupy upon arrival the server is different compared with the service times of the retrial jobs. Moreover, the service requirements of each class of retrial customers is also different. Besides its practical applicability in the modelling of relay assisted cooperative wireless networks, and in call centers with call-back option, our work is also theoretically oriented. 

In particular, in this work we focus on the fundamental
problem of investigating the queueing delay in multiclass retrial systems with constant retrial policy, and arbitrarily distributed class dependent service times, which remains an open problem. The only available results refer to the investigation of the stability conditions \cite{avr3,avr5,mor,mord}. More precisely, for the two orbit scenario, we generalize the seminal paper in \cite{avr} by allowing arbitrarily distributed class dependent service times, and obtain the generating function of the stationary joint orbit queue-length distribution in terms of a solution of a Riemann boundary value problem\footnote{In subsection \ref{fred} we also provided the way we can expressed it by solving a Fredholm integral equation of the second kind.}. Our contribution provides a building block towards the generalization to the case of $N$ orbits; see also Section \ref{conc}. For the completely symmetrical system, we also provide for the first time, explicit expressions for the expected number of customers at each orbit queue, without the need of solving a boundary value problem.

The rest of the paper is organized as follows. In Section \ref{sec:1} we describe the model in detail and provide the fundamental functional equation. Some important preparatory results are given in Section \ref{gener}. Sections \ref{msym}, and \ref{asym} are devoted in the detailed analysis of the modified symmetrical and the asymmetrical system, respectively. In Section \ref{sym} we provide explicit expressions for the expected orbit delay for the completely symmetrical system without solving a boundary value problem, while in Section \ref{num} a simple numerical example is presented.
\section{The model}\label{sec:1}
Consider a single server queue accepting two types of customers, say $P_{1}$, $P_{2}$. $P_{i}$, customers arrive according to Poisson process with rate $\lambda_{i}$, and if upon arrival find the server unavailable, enter a dedicated virtual queue, called the orbit queue $i$, $i=1,2$. All the customers in each orbit behave independently of each other and try to access the server according to the constant retrial policy. More precisely, we assume that the retrial times for any orbiting $P_{i}$ customer are
exponentially distributed with rate $\theta_{i}/n_{i}$, given that there are $n_{i}$ customers in orbit $i$, $i=1,2,$. Upon a service completion, the server remains idle until either a primary or a retrial customer (of either type) arrives.

The provided service time depends on the type (i.e., $P_{1}$, $P_{2}$) and the state of the customer (i.e., either orbiting or primary). More precisely, the service times for orbiting customers of type $i$, say $B_{i},\,i=1,2,$ is arbitrarily distributed with cumulative distribution function (cdf) $B_{i}(x)$, probability density function (pdf) $b_{i}(x)$, Laplace Stieltjes Transform (LST) $\beta_{i}^{*}(s)$, and moments $\bar{b}_{i}$, $\bar{b}_{i}^{(2)}$. An arriving primary customer of either type who finds the server idle will occupy it immediately and its service requirement, say $B_{3}$, is arbitrarily distributed with cdf $B_{3}(x)$, pdf $b_{3}(x)$, LST $\beta_{3}^{*}(s)$, and moments $\bar{b}_{3}$, $\bar{b}_{3}^{(2)}$. 

Let $X_{i}(n)$ be the number of $P_{i}$, $i=1,2,$ orbiting customers, just after the end of the $n$th service completion. Denote also by $\xi(n)$, the type of the $n$th service time. Clearly $X(n)=(X_{1}(n), X_{2}(n), \xi(n))$ forms an irreducible and aperiodic Markov chain. Define by $A_{ij}(n)$, $i,j=1,2,3$ the number of $P_{i}$ customers that arrive during the $n$th service service period if it is of type $j$. Then,
\begin{displaymath}
\begin{array}{l}
(X_{1}(n+1),X_{2}(n+1),\xi(n+1))\\=\begin{cases}
(X_{1}(n)-1+A_{11}(n+1), X_{2}(n)+A_{21}(n+1), 1), & w.p. \frac{\theta_{1}}{D(n)}\\
(X_{1}(n)+A_{13}(n+1), X_{2}(n)+A_{23}(n+1), 3), & w.p. \frac{\lambda}{D(n)}, \\
(X_{1}(n)+A_{12}(n+1), X_{2}(n)-1+A_{22}(n+1),2), & w.p. \frac{\theta_{2}}{D(n)},
\end{cases}
\end{array}
\end{displaymath}
where $D(n)=\lambda+\theta_{1}\textbf{1}_{X_{1}(n)>0}+\theta_{2}\textbf{1}_{X_{2}(n)>0}$ and $\lambda=\lambda_{1}+\lambda_{1}$. Denote,
\begin{displaymath}
\begin{array}{rl}
\pi_{m,l}=&\lim_{n\rightarrow \infty}Pr((X_{1}(n),X_{2}(n))=(m,l)),\\ \Pi(z_{1},z_{2})=&\sum_{m=0}^{\infty}\sum_{l=0}^{\infty}\pi_{m,l}z_{1}^{m}z_{2}^{l},\,|z_{1}|\leq1, |z_{2}|\leq1,
\end{array}
\end{displaymath}
and $\theta=\theta_{1}+\theta_{2}$, $r_{i}=\lambda_{i}/\lambda$. Clearly, for $j=1,2,3,$
\begin{displaymath}
\begin{array}{rl}
P(A_{1j}=k,A_{2j}=m)=d_{km}^{(j)}=&\int_{0}^{\infty}e^{-\lambda_{1}t}\frac{(\lambda_{1}t)^{k}}{k!}e^{-\lambda_{2}t}\frac{(\lambda_{2}t)^{k}}{m!}dB_{j}(x),\\
d_{j}^{*}(z_{1},z_{2})=\sum_{m=0}^{\infty}\sum_{l=0}^{\infty}d_{km}^{(j)}z_{1}^{k}z_{2}^{m}=&\beta_{j}^{*}(\lambda(1-r_{1}z_{1}-r_{2}z_{2})).
\end{array}
\end{displaymath}
Let $y=\lambda(1-r_{1}z_{1}-r_{2}z_{2})$. Considering the transition probabilities at service completion epochs we obtain,
\begin{displaymath}
\begin{array}{rl}
\pi_{m,l}=&\frac{\theta_{1}}{\lambda+\theta}\sum_{k_{1}=1}^{m+1}\sum_{k_{2}=1}^{l}\pi_{k_{1},k_{2}}d_{m+1-k_{1},l-k_{2}}^{(1)}+\frac{\theta_{2}}{\lambda+\theta}\sum_{k_{1}=1}^{m}\sum_{k_{2}=1}^{l+1}\pi_{k_{1},k_{2}}d_{m-k_{1},l+1-k_{2}}^{(2)}\vspace{2mm}\\
&\frac{\lambda}{\lambda+\theta}\sum_{k_{1}=1}^{m}\sum_{k_{2}=1}^{l}\pi_{k_{1},k_{2}}d_{m-k_{1},l-k_{2}}^{(3)}+\frac{\theta_{1}}{\lambda+\theta_{1}}\sum_{k_{1}=1}^{m+1}\pi_{k_{1},0}d_{m+1-k_{1},l}^{(1)}+\pi_{0,0}d_{m,l}\vspace{2mm}\\
&+\frac{\lambda}{\lambda+\theta_{1}}\sum_{k_{1}=1}^{m}\pi_{k_{1},0}d_{m-k_{1},l}^{(3)}+\frac{\theta_{2}}{\lambda+\theta_{2}}\sum_{k_{2}=1}^{l+1}\pi_{0,k_{2}}d_{m,l+1-k_{2}}^{(2)}+\frac{\lambda}{\lambda+\theta_{2}}\sum_{k_{2}=1}^{l}\pi_{0,k_{2}}d_{m,l-k_{2}}^{(2)}.
\end{array}
\end{displaymath}
Forming the generating functions we conclude that
\begin{equation}
K(z_{1},z_{2})\Pi(z_{1},z_{2})=A(z_{1},z_{2})\Pi(z_{1},0)+B(z_{1},z_{2})\Pi(0,z_{2})+C(z_{1},z_{2})\Pi(0,0),
\label{equ}
\end{equation}
where,
\begin{equation}
K(z_{1},z_{2})=z_{1}z_{2}-\widetilde{K}(z_{1},z_{2}).
\label{ker}
\end{equation}
\begin{displaymath}
\begin{array}{rl}
A(z_{1},z_{2})=&z_{2}\widetilde{A}(z_{1},z_{2})-\widetilde{K}(z_{1},z_{2}),\\
B(z_{1},z_{2})=&z_{1}\widetilde{B}(z_{1},z_{2})-\widetilde{K}(z_{1},z_{2}),\\
C(z_{1},z_{2})=&\widetilde{K}(z_{1},z_{2})+z_{2}(r_{1}z_{1}\beta_{3}^{*}(y)-\widetilde{A}(z_{1},z_{2}))+z_{1}(r_{2}z_{2}\beta_{3}^{*}(y)-\widetilde{B}(z_{1},z_{2})),
\end{array}
\end{displaymath}
and
\begin{displaymath}
\begin{array}{rl}
\widetilde{K}(z_{1},z_{2})=&\frac{\theta_{1}}{\lambda+\theta}z_{2}\beta_{1}^{*}(y)+\frac{\theta_{2}}{\lambda+\theta}z_{1}\beta_{2}^{*}(y)+\frac{\lambda z_{1}z_{2}}{\lambda+\theta}\beta_{3}^{*}(y),\\
\widetilde{A}(z_{1},z_{2})=\frac{\theta_{1}\beta_{1}^{*}(y)+\lambda z_{1}\beta_{3}^{*}(y)}{\lambda+\theta_{1}},&
\widetilde{B}(z_{1},z_{2})=\frac{\theta_{2}\beta_{2}^{*}(y)+\lambda z_{2}\beta_{3}^{*}(y)}{\lambda+\theta_{2}}.
\end{array}
\end{displaymath}
$K(z_{1},z_{2})$ is called the kernel of the functional equation (\ref{equ}), and its investigation is of major importance for the fruitful analysis of (\ref{equ}). Contrary to \cite{box1,dim}, $K(z_{1},z_{2})$ is not a Poisson kernel.
\section{General results}\label{gener}
Some interesting results can be deduced directly by the functional equation. Substituting $z_{1}=1$ in (\ref{equ}) and subsequently letting $z_{2}\rightarrow 1$, and vice versa yield the following linear relations between $\Pi(0,1)$, $\Pi(1,0)$ and $\Pi(0,0)$.
\begin{equation}
\begin{array}{rl}
1-\widehat{\rho}_{2}=&\Pi(1,0)\frac{\lambda+\theta_{1}+\lambda_{2}[\lambda(\overline{b}_{3}-\overline{b}_{2})+\theta_{1}(\overline{b}_{1}-\overline{b}_{2})]}{\lambda+\theta_{1}}+\Pi(0,1)\frac{\theta_{1}}{\theta_{2}}[\frac{\lambda_{2}[\lambda(\overline{b}_{3}-\overline{b}_{1})+\theta_{2}(\overline{b}_{2}-\overline{b}_{1})]-\theta_{2}}{\lambda+\theta_{2}}]\\
&+\Pi(0,0)\theta_{1}[\frac{\lambda_{2}(\overline{b}_{3}-\overline{b}_{1})}{\lambda+\theta_{1}}+\frac{1+\lambda_{2}(\overline{b}_{3}-\overline{b}_{2})}{\lambda+\theta_{2}}],\vspace{2mm}\\
1-\widehat{\rho}_{1}=&\Pi(1,0)\frac{\theta_{2}}{\theta_{1}}[\frac{\lambda_{1}[\lambda(\overline{b}_{3}-\overline{b}_{2})+\theta_{1}(\overline{b}_{1}-\overline{b}_{2})]-\theta_{1}}{\lambda+\theta_{1}}]+\Pi(0,1)\frac{\lambda+\theta_{2}+\lambda_{1}[\lambda(\overline{b}_{3}-\overline{b}_{1})+\theta_{2}(\overline{b}_{2}-\overline{b}_{1})]}{\lambda+\theta_{2}}\\
&+\Pi(0,0)\theta_{2}[\frac{\lambda_{1}(\overline{b}_{3}-\overline{b}_{2})}{\lambda+\theta_{2}}+\frac{1+\lambda_{1}(\overline{b}_{3}-\overline{b}_{1})}{\lambda+\theta_{1}}],
\end{array}
\label{qa}
\end{equation}
where $\widehat{\rho}_{j}=\frac{\lambda_{j}(\theta_{1}\overline{b}_{1}+\theta_{2}\overline{b}_{2}+\lambda\overline{b}_{3})}{\theta_{j}}$, $j=1,2.$

We proceed with an interesting interpretation for $\widehat{\rho}_{j}$. Let $S_{j}$, $j=1,2$ be the time elapsed form the epoch a service is initiated until the epoch the server becomes idle after a service completion of a retrial customer of type $j$ given that both orbit queues are non-empty, and $N_{i}(S_{j})$ the number of type $i$ customers that join the orbit queue $i$ during $S_{j}$. Let also $s^{(j)}_{k_{1},k_{2}}(t)dt=P(t<S_{j}\leq t+dt,N_{i}(S_{j})=k_{i})$. We restrict the analysis to the orbit queue 1. The analysis for the orbit queue 2 is similar. Then,
\begin{displaymath}
\begin{array}{rl}s^{(1)}_{k_{1},k_{2}}(t)=&\frac{\theta_{1}}{\lambda+\theta}h_{k_{1},k_{2}}(t)b_{1}(t)+\frac{\lambda}{\lambda+\theta}\sum_{m_{1}=0}^{k_{1}}\sum_{m_{2}=0}^{k_{2}}h_{k_{1},k_{2}}(t)b_{3}(t)*s^{(j)}_{k_{1}-m_{1},k_{2}-m_{2}}(t)\\
&+\frac{\theta_{2}}{\lambda+\theta}\sum_{m_{1}=0}^{k_{1}}\sum_{m_{2}=0}^{k_{2}}h_{k_{1},k_{2}}(t)b_{2}(t)*s^{(j)}_{k_{1}-m_{1},k_{2}-m_{2}}(t),
\end{array}
\end{displaymath}
where $h_{m,n}(t)=e^{-\lambda t}\frac{(\lambda_{1}t)^{m}}{m!}\frac{(\lambda_{2}t)^{n}}{n!}$ and ``*'' means convolution. If 
\begin{displaymath}
\begin{array}{c}
\widetilde{s}_{j}^{*}(z_{1},z_{2},s)=\int_{0}^{\infty}e^{-st}\sum_{k_{1}=0}^{\infty}\sum_{k_{2}=0}^{\infty}s_{k_{1},k_{2}}^{(j)}(t)z_{1}^{k_{1}}z_{2}^{k_{2}}dt,
\end{array}
\end{displaymath}
then,
\begin{displaymath}
\begin{array}{rl}\widetilde{s}_{1}^{*}(z_{1},z_{2},s)=&\frac{\theta_{1}\beta_{1}^{*}(s+y)}{s+\theta_{1}+\theta_{2}(1-\beta_{2}^{*}(s+y))+\lambda(1-\beta_{3}^{*}(s+y))},\\
\widetilde{s}_{2}^{*}(z_{1},z_{2},s)=&\frac{\theta_{2}\beta_{2}^{*}(s+y)}{s+\theta_{2}+\theta_{1}(1-\beta_{2}^{*}(s+y))+\lambda(1-\beta_{3}^{*}(s+y))},
\end{array}
\end{displaymath}
and
\begin{displaymath}
\widehat{\rho}_{1}=\frac{\partial}{\partial z_{1}}\widetilde{s}_{1}^{*}(z_{1},1,0)|_{z_{1}=1},\,\widehat{\rho}_{2}=\frac{\partial}{\partial z_{2}}\widetilde{s}_{2}^{*}(1,z_{2},0)|_{z_{2}=1}.
\end{displaymath}
That said, $\widehat{\rho}_{j}$ is the expected number of customers that join the orbit queue $j$ during this special service time $S_{j}$. Therefore, we expect that $\widehat{\rho}_{j}<1$, $j=1,2$, which is consistent with the results regarding the stability conditions derived in \cite{avr3}. 
\subsection{Special cases}
\paragraph{The modified symmetrical model}Consider the modified symmetrical model where, $\lambda_{1}=\lambda_{2}=\frac{\lambda}{2}$ (i.e., $r_{1}=r_{2}$), $\theta_{1}=\theta_{2}=\frac{\theta}{2}$ and $B_{1}\sim B_{2}\sim B$ and $B_{3}\nsim B$. Then, (\ref{qa}) becomes
\begin{equation}
\begin{array}{rl}
1-\frac{\lambda(\theta\overline{b}+\lambda\overline{b}_{3})}{\theta}=&\Pi(1,0)\frac{\theta+\lambda(2+\lambda(\overline{b}_{3}-\overline{b}))}{2\lambda+\theta}
+\Pi(0,1)\frac{\lambda^{2}(\overline{b}_{3}-\overline{b})-\theta}{2\lambda+\theta}\\
&+\Pi(0,0)\frac{\theta(1+\lambda(\overline{b}_{3}-\overline{b}))}{2\lambda+\theta}\\
1-\frac{\lambda(\theta\overline{b}+\lambda\overline{b}_{3})}{\theta}=&\Pi(1,0)\frac{\lambda^{2}(\overline{b}_{3}-\overline{b})-\theta}{2\lambda+\theta}
+\Pi(0,1)\frac{\theta+\lambda(2+\lambda(\overline{b}_{3}-\overline{b}))}{2\lambda+\theta}\\
&+\Pi(0,0)\frac{\theta(1+\lambda(\overline{b}_{3}-\overline{b}))}{2\lambda+\theta}.\end{array}
\label{qaa}
\end{equation}
By subtracting the above equations we conclude that $\Pi(0,1)=\Pi(1,0)$ and substituting back we derive,
\begin{equation}
\begin{array}{c}
\frac{(1-\widehat{\rho})(2\lambda+\theta)}{1+\lambda(\overline{b}_{3}-\overline{b})}=2\lambda\Pi(1,0)+\theta\Pi(0,0).
\end{array}
\label{df1}
\end{equation}
Since the right hand side of the above equation is positive, it is straightforward that $\widehat{\rho}=\frac{\lambda(\theta\overline{b}+\lambda\overline{b}_{3})}{\theta}<1$ is the ergodicity condition.

\paragraph{The completely symmetrical model}Let $\lambda_{1}=\lambda_{2}=\frac{\lambda}{2}$ (i.e., $r_{1}=r_{2}$), $\theta_{1}=\theta_{2}=\frac{\theta}{2}$, $B_{j}\sim B$, $j=1,2,3$. Equations (\ref{qaa}), (\ref{df1}) remain valid with slight modifications, and the stability condition is $\widehat{\rho}=\frac{\lambda(\theta+\lambda)\overline{b}}{\theta}<1\Leftrightarrow\theta-\lambda\bar{b}(\lambda+\theta)>0.$
\section{Detailed analysis of the modified symmetrical model}\label{msym}
Consider the modified symmetrical model, where $\lambda_{1}=\lambda_{2}=\frac{\lambda}{2}$ (i.e., $r_{1}=r_{2}=\frac{1}{2}$), $\theta_{1}=\theta_{2}=\frac{\theta}{2}$, $B_{1}\sim B_{2}\sim B$ and $B_{3}\nsim B$, and assume that $\widehat{\rho}=\frac{\lambda}{\theta}(\theta\overline{b}+\lambda\overline{b}_{3})<1$.
\subsection{Preliminary analysis}
We now follow the methodology given in \cite{bv}. Let,
\begin{displaymath}
\begin{array}{c}w_{i}=2r_{i}z_{i}=z_{i},\,i=1,2,\,\delta=\frac{1}{2}(w_{1}+w_{2})=\frac{1}{2}(z_{1}+z_{2}).\end{array}
\end{displaymath}
Clearly,
\begin{displaymath}
\begin{array}{c}
K(z_{1},z_{2})=z_{1}z_{2}-\frac{\theta\beta^{*}(y)}{2(\lambda+\theta)}(z_{1}+z_{2})-\frac{\lambda z_{1}z_{2}}{\lambda+\theta}\beta_{3}^{*}(y),
\end{array}
\end{displaymath}
is well defined for $z_{1},z_{2}$ with $Re(\delta)\leq 1$ and $K(z_{1},2\delta-z_{1})=0$ is a quadratic equation in $z_{1}$ for every fixed $\delta$ with $\Re(\delta)\leq 1$. Particularly,
\begin{displaymath}
\begin{array}{c}
K(z_{1},2\delta-z_{1})=0\Rightarrow z_{1}^{2}-2\delta z_{1}+\delta\check{\beta}^{*}(\delta)=0,
\end{array}
\end{displaymath}
where $\check{\beta}^{*}(\delta)=\frac{\theta\beta^{*}(\lambda(1-\delta))}{\theta+\lambda(1-\beta_{3}^{*}(\lambda(1-\delta)))}$,
and it has two roots $\widehat{z}_{1}=\widehat{z}_{1}(\delta)$, $\widehat{z}_{2}=\widehat{z}_{1}(\delta)=2\delta-\widehat{z}_{1}(\delta)$. The equation $K(z_{1},2\delta-z_{1}) = 0$ can also be written as
\begin{equation}
\begin{array}{c}
(z_{1}-\delta)^{2}=(\delta-\frac{\check{\beta}^{*}(\delta)}{2})^{2}-(\frac{\check{\beta}^{*}(\delta)}{2})^{2}.
\end{array}
\label{fgh}
\end{equation}
It is easy to check that the right hand side of (\ref{fgh}) is the determinant of $K(z_{1},2\delta-z_{1}) = 0$, given by $D(\delta)=\delta(\delta-\check{\beta}^{*}(\delta))$. Clearly, $D(\delta)=0$ has two roots in $\Re(\delta)\leq 1$, viz. $\delta_{0}=0$ and $\delta_{1}=1$, since $\delta-\check{\beta}^{*}(\delta)$ has exactly one zero (i.e., $\delta_{1}$) in $\Re(\delta)\leq 1$ when $\widehat{\rho}=\lambda\frac{\theta\overline{b}+\lambda\overline{b}_{3}}{\theta}<1$.

Note now that $\check{\beta}^{*}(z_{1},z_{2})=\frac{\theta\beta^{*}(y)}{\theta+\lambda(1-\beta_{3}^{*}(y))}$ has a very intuitive probabilistic interpretation. Indeed, let $\check{B}$ be the time elapsed from the epoch a service is initiated, until the service completion of a retrial customer of either type, given that both orbit queues are non-empty. Let $N_{i}(\check{B})$ is the number of newly arriving type $i$ customers during $\check{B}$. Then,
\begin{equation}
\begin{array}{rl}
\check{b}_{l_{1},l_{2}}(t)dt=&P(t<\check{B}\leq t+dt,N_{i}(\check{B})=l_{i},i=1,2),\\
\check{b}_{l_{1},l_{2}}(t)=&\frac{\theta}{\lambda+\theta}h_{l_{1},l_{2}}(t)b(t)+\frac{\lambda}{\lambda+\theta}\sum_{m_{1}=0}^{l_{1}}\sum_{m_{2}=0}^{l_{2}}h_{m_{1},m_{2}}(t)
b(t)*\check{b}_{l_{1}-m_{1},l_{2}-m_{2}}(t),
\end{array}
\label{1fgh}
\end{equation}
where ``$*$'' means convolution. If $\check{\beta}^{*}(z_{1},z_{2},s)=\int_{0}^{\infty}\sum_{l_{1}=0}^{\infty}\sum_{l_{2}=0}^{\infty}\check{b}_{l_{1},l_{2}}(t)z_{1}^{l_{1}}z_{2}^{l_{2}}$, $|z_{i}|<1$, $i=1,2,$ we have $\check{\beta}(z_{1},z_{2})=\check{\beta}^{*}(z_{1},z_{2},0)$.

Put $G=\left\{\delta:0\leq\delta\leq1\right\}$, and consider the two-bladed Riemann surface $S$ composed of two semi-planes $\left\{\delta:\Re(\delta)\leq 1\right\}$ slitted along $G$, then $\widehat{z}_{1}(\delta)$ and also $\widehat{z}_{2}(\delta)$ constitute analytic functions on $S$ for $\Re(\delta)< 1$.

Next we introduce the following parametrization of $G$. Consider the function
\begin{equation}
\begin{array}{c}
\delta-\frac{\check{\beta}^{*}(\delta)}{2}(1+\cos\phi),\,\Re(\delta)\leq1,\,\phi\in[0,2\pi].
\end{array}
\label{poo}
\end{equation}
Using  Rouche's theorem it can be proved that the function in (\ref{poo}) has exactly one zero, say $\delta=\delta(\phi)$ in $\Re(\delta)\leq 1$ for $\phi\in[0,2\pi]$, which is real. Thus, $G=\left\{\delta:\delta=\delta(\phi),\,\phi\in[0,2\pi]\right\}$. Therefore, for $\phi\in[0,2\pi]$, substitute the zero $\delta=\delta(\phi)$ of (\ref{poo}) in (\ref{fgh}), we have
\begin{equation}
\begin{array}{rl}
\widehat{z}_{1}=\widehat{z}_{1}(\delta(\phi))=&\frac{\check{\beta}^{*}(\delta(\phi))}{2}(1+e^{i\phi}),\\
\widehat{z}_{2}=\widehat{z}_{2}(\delta(\phi))=&\frac{\check{\beta}^{*}(\delta(\phi))}{2}(1+e^{-i\phi}).
\end{array}
\label{ww1}
\end{equation}
Let $L_{1}=\left\{z_{1}: z_{1}=\widehat{z}_{1}(\delta(\phi));\ \phi\in[0,2\pi]\right\},\,L_{2}=\left\{z_{2}: z_{2}=\widehat{z}_{2}(\delta(\phi));\ \phi\in[0,2\pi]\right\}$. Then, the following statements are readily verified:
$i)$ $L_{1}$, and similarly $L_{2}$, is a simple smooth contour, $ii)$ $L_{1}\subset \left\{z_{1}:|z_{1}|\leq 1\right\}$, $L_{2}\subset \left\{z_{2}:|z_{2}|\leq 1\right\}$, $iii)$ $z_{1}=1\in L_{1}$, $z_{1}=0\in L_{1}^{+}$, $z_{2}=1\in L_{2}$, $z_{2}=0\in L_{2}^{+}$, $iv)$ The relations in (\ref{ww1}) define a one to one mapping $\widehat{z}_{1}(\delta(\phi))=\omega_{1}(\widehat{z}_{2}(\delta(\phi)))$ of $L_{2}$ onto $L_{1}$.

Clearly, the contours $L_{1}$, $L_{2}$ satisfy the conditions of Theorem 1.1 in \cite{coh2}, p. 101. Put, for $\phi\in[0,2\pi]$,
\begin{displaymath}
\begin{array}{c}\rho(\phi)=|\widehat{z}_{1}(\delta(\phi))|,\,\omega(\phi)=\arg(\widehat{z}_{1}(\delta(\phi)),\end{array}
\end{displaymath}
so we can write: $\widehat{z}_{1}(\delta(\phi))=\rho(\phi)e^{i\omega(\phi)},\,\widehat{z}_{2}(\delta(\phi))=\rho(\phi)e^{-i\omega(\phi)}.$

We proceed by applying Theorem 1.1 in \cite{coh2}, and thus there exists a unique simple contour $L$ in the $z-$plane with
$z=0\in L^{+},\,z=1\in L,\,z=\infty\in L^{-}$,
and functions
$f_{1}(z):L^{+}\cup L\rightarrow L_{1}^{+}\cup L_{1},\,f_{2}(z):L^{-}\cup L\rightarrow L_{2}^{+}\cup L_{2}$,
such that: $i)$ $z=0$ is a simple zero of $f_{1}(.)$, $ii)$ $z=\infty$ is a simple zero of $f_{2}(.)$, i.e., $0<d:=\,\lim_{|z| \to \infty}|zf_{2}(z)|<\infty$, $iii)$ $f_{1}:L^{+}\rightarrow L_{1}^{+}$ is regular and univalent for $z\in L^{+}$, $iv)$ $f_{2}:L^{-}\rightarrow L_{2}^{+}$ is regular and univalent for $z\in L^{-}$, $v)$ $f_{1}(z)=\omega_{1}(f_{2}(z))$, $z\in L$.

Therefore, $\log (\frac{f_{1}(z)}{z})$ should be regular for $z\in L^{+}$ and continuous for $z\in L^{+}\cup L$, and $\log(zf_{2}(z))$ should be regular for $z\in L^{-}$ and continuous for $z\in L^{-}\cup L$. Let $\psi(z)$, $z\in L$ with $\psi(1)=0$, be a real function with $\phi=\psi(z)$. Then,
\begin{displaymath}
\begin{array}{c}
(\frac{f_{1}(z)}{z})+\frac{\log(zf_{2}(z))}{d}=\frac{\log(\rho^{2}(\psi(z)))}{d},\,z\in L.
\end{array}
\end{displaymath}

If $\log(\rho^{2}(\psi(z)))$ satisfies the Holder condition on $L$, the equation above represent a simple Riemann boundary value problem and following \cite{bv}, \cite{coh2},
\begin{equation}
\begin{array}{rl}
f_{1}(z)=&z\exp\left\{\frac{1}{2i\pi}\int_{\zeta\in L}\log(\rho(\psi(\zeta)))[\frac{\zeta+z}{\zeta-z}-\frac{\zeta+1}{\zeta-1}]\frac{d\zeta}{\zeta}\right\},\,z\in L^{+},\\
f_{2}(z)=&z^{-1}\exp\left\{-\frac{1}{2i\pi}\int_{\zeta\in L}\log(\rho(\psi(\zeta)))[\frac{\zeta+z}{\zeta-z}-\frac{\zeta+1}{\zeta-1}]\frac{d\zeta}{\zeta}\right\},\,z\in L^{-}.
\end{array}
\label{ww2}
\end{equation}
By applying the Plemelj-Sokhotski formulas we obtain,
\begin{equation}
\begin{array}{rl}
f_{1}(z)=&z\rho(\psi(z))\exp\left\{\frac{1}{2i\pi}\int_{\zeta\in L}\log(\rho(\psi(\zeta)))[\frac{\zeta+z}{\zeta-z}-\frac{\zeta+1}{\zeta-1}]\frac{d\zeta}{\zeta}\right\},\,z\in L,\\
f_{2}(z)=&z^{-1}\rho(\psi(z))\exp\left\{-\frac{1}{2i\pi}\int_{\zeta\in L}\log(\rho(\psi(\zeta)))[\frac{\zeta+z}{\zeta-z}-\frac{\zeta+1}{\zeta-1}]\frac{d\zeta}{\zeta}\right\},\,z\in L.
\end{array}
\label{ww3}
\end{equation}
From these expressions it is seen (see \cite{bv}, p. 99) that $\psi(z)$ should satisfy,
\begin{equation}
\begin{array}{l}
e^{i\omega(\psi(z))}=z\exp\left\{\frac{1}{2i\pi}\int_{\zeta\in L}\log(\rho(\psi(\zeta)))[\frac{\zeta+z}{\zeta-z}-\frac{\zeta+1}{\zeta-1}]\frac{d\zeta}{\zeta}\right\},\,z\in L.
\end{array}
\label{lp}
\end{equation}
The solution of the Riemann problem above depends on the value of the constant $d$, which is chosen such that $z = 1\in L$. Thus, $f_{1}(1)$ corresponds $f_{1}(\delta(0))$, and $f_{2}(1)$ to $f_{2}(\delta(0))$.
 
We proceed with the solution of the functional equation. Since $(\widehat{z}_{1},\widehat{z}_{2})$ with $\widehat{z}_{1}=f_{1}(z)$, $\widehat{z}_{2}=f_{2}(z)$, $z\in L\equiv\left\{z:|z|=1\right\}$\footnote{This is due to the symmetry of the model.} (see Theorem 4.1 in \cite{bv} or Section 4 in\cite{coh2}), is a zerotuple of the kernel $K(z_{1}, z_{2})$ with $f_{1}(z)$, $z\in L^{+}\cup L$, $f_{2}(z)$, $z \in L^{-}\cup L$ as constructed above, it follows that for $|z|=1$,
\begin{equation}
\begin{array}{c}
\Pi(f_{1}(z),0)A(f_{1}(z),f_{2}(z))+\Pi(0,f_{2}(z))B(f_{1}(z),f_{2}(z))+\Pi(0,0)C(f_{1}(z),f_{2}(z))=0.
\end{array}
\label{rui}
\end{equation}
Moreover, it follows from the regularity of $f_{1}(z)$, $z\in L^{+}$, $f_{2}(z)$, $z\in L^{-}$, that $i)$ $\widehat{\Pi}_{1}(z)=\Pi(f_{1}(z),0)/\Pi(0,0)$, $z\in L^{+}\cup L$ is regular for $z\in L^{+}$ and continuous for $L^{+}\cup L$, $ii)$ $\widehat{\Pi}_{2}(z)=\Pi(0,f_{2}(z))/\Pi(0,0)$, $z\in L^{-}\cup L$ is regular for $z\in L^{-}$ and continuous for $L^{-}\cup L$, $iii)$ $\widehat{\Pi}_{1}(0)=1$, $\lim_{|z|\rightarrow \infty}\widehat{\Pi}_{2}(z)=1$.

Note also that $|f_{1}(z)|\leq 1$ for $|z|=1$, so that the regularity of $f_{1}(z)$ for $z\in L^{+}$ implies by means of the maximum modulus theorem that $|f_{1}(z)|< 1$ for $|z|<1$, so that $\widehat{\Pi}_{1}(z)$ is well defined, analogously for $\widehat{\Pi}_{2}(z)$. Furthermore, $\widehat{\Pi}_{1}(0)=1$, $\lim_{|z|\rightarrow \infty}\widehat{\Pi}_{2}(z)=1$. Then, (\ref{rui}) can be rewritten as
\begin{equation}
\begin{array}{c}
\widehat{\Pi}_{1}(z)=G(z)\widehat{\Pi}_{2}(z)+g(z),\,|z|=1,
\end{array}
\label{roi}
\end{equation}
where now,
\begin{equation}
\begin{array}{rl}
G(z)=&-\frac{f_{1}(z)}{f_{2}(z)}\frac{f_{2}(z)-\widehat{\beta}(f_{1}(z),f_{2}(z))}{f_{1}(z)-\widehat{\beta}(f_{1}(z),f_{2}(z))},\\
g(z)=&-G(z)+1+\frac{(\lambda+\frac{\theta}{2})f_{1}(z)(1-\beta_{3}^{*}(\lambda(1-\frac{f_{1}(z)+f_{2}(z)}{2})))}{(f_{1}(z)-\widehat{\beta}(f_{1}(z),f_{2}(z)))(\frac{\theta}{2}
+\lambda(1-\beta_{3}^{*}(\lambda(1-\frac{f_{1}(z)+f_{2}(z)}{2}))))},
\end{array}
\label{rui1}
\end{equation}
where $\widehat{\beta}(z_{1},z_{2})=\frac{\theta\beta^{*}(y)}{2(\frac{\theta}{2}+\lambda(1-\beta_{3}^{*}(y)))}$. Using similar arguments as in the derivation of (\ref{1fgh}) we can easily prove that $\widehat{\beta}(z_{1},z_{2})$ is a probability generating function of a proper probability mass function.

Clearly, (\ref{roi}) along with the above conditions to be satisfied by $\widehat{\Pi}_{1}$, $\widehat{\Pi}_{2}$ formulate a Riemann boundary value problem. For its analysis we have firstly to discuss some properties of $G$ and $g$. From the definition of $f_{1}(z)$, $f_{2}(z)$ we have,
\begin{displaymath}
\begin{array}{l}
f_{1}(z)=\widehat{z}_{1}(\delta(\psi(z))),\,f_{2}(z)=\widehat{z}_{2}(\delta(\psi(z))),\,
0<f_{1}(z)+f_{2}(z)\leq2,\,|z|=1.
\end{array}
\end{displaymath}
Consequently, $0<\widehat{\beta}(f_{1}(z),f_{2}(z))\leq1.$
Furthermore, it is not difficult to show that $1-f_{1}(z)$, $1-f_{2}(z)$ have a zero of multiplicity one at $z=1$, and also that $1-\beta_{3}^{*}(\lambda(1-\frac{f_{1}(z)+f_{2}(z)}{2}))$ and $f_{1}(z)-\widehat{\beta}(f_{1}(z),f_{2}(z))$ have a zero at $z=1$. Thus, it follows that $G(1)$, $g(1)$ are bounded (we have to note here that $z=1\in L$ and $f_{i}(1)$ corresponds to $f_{i}(\delta(0))$, $i=1,2$. Thus, $G(1)\equiv G(\delta(\psi(1)))$).

Moreover, the other point of interest is $\delta(\pi)$. Clearly $\widehat{z}_{1}(\delta(\pi))=0=\widehat{z}_{2}(\delta(\pi))$, and as a result the numerator and the denominator of $G(\delta(\pi))$ vanish simultaneously. Thus, $\delta(\pi)$ is a cancelled point of $G(\delta(\phi))$. Therefore, $0<G(\delta(\phi))<\infty$. Similarly we can prove that $0<g(\delta(\phi))<\infty$. To conclude $G(z)$, $g(z)$ never vanishes for $|z|=1$.

Clearly, we can easily show starting by (\ref{ww3}) that $G(z)$ and also $g(z)$ both possess a continuous derivative along $|z| = 1$ (note that $L_{1}$, $L_{2}$ and $L$ are all smooth contours) and consequently, they satisfy the Holder condition on $|z| =1$.
\subsection{Solution of a Riemann boundary value problem}\label{rh}
In order to solve the Riemann boundary value problem formulated by (\ref{roi}) and conditions 1, 2, we have to compute the index of $G(z)$ on $|z|=1$. Note that
\begin{displaymath}
\begin{array}{l}
\chi=ind_{|z|=1}G(z)=ind_{|z|=1}f_{1}(z)-ind_{|z|=1}f_{2}(z)\\+ind_{|z|=1}[f_{2}(z)-\widehat{\beta}(f_{1}(z),f_{2}(z))]
-ind_{|z|=1}[f_{1}(z)-\widehat{\beta}(f_{1}(z),f_{2}(z))].
\end{array}
\end{displaymath}

Since $L_{1}$ and $L_{2}$ are simple contours with $z_{1}=0 \in L_{1}^{+}$, $z_{2} = 0\in L_{2}^{+}$, and $f_{1}(z)$ traverses $L_{1}$ counterclockwise, whereas $f_{2}(z)$ traverses $L_{2}$ clockwise we have $ind_{|z|=1}f_{1}(z)=1$, $ind_{|z|=1}f_{2}(z)=-1$. The contours $f_{j}(z)-\widehat{\beta}(f_{1}(z),f_{2}(z))$, are smooth and have only two real points of which one is negative and the other that corresponds to $z=1$ is located at zero where the contours have vertical tangents. Thus,
\begin{displaymath}
\begin{array}{rl}
ind_{|z|=1}[f_{2}(z)-\widehat{\beta}(f_{1}(z),f_{2}(z))]=&-\frac{1}{2},\\
ind_{|z|=1}[f_{1}(z)-\widehat{\beta}(f_{1}(z),f_{2}(z))]=&\frac{1}{2}.
\end{array}
\end{displaymath}
Therefore, $\chi=1$ and,
\begin{equation}
\begin{array}{rl}
\widehat{\Pi}_{1}(z)=&e^{\Gamma_{1}(z)}[\Psi(z)+c_{1}z+c_{0}],\,|z|<1,\\
\widehat{\Pi}_{2}(z)=&z^{-1}e^{\Gamma_{1}(z)}[\Psi(z)+c_{1}z+c_{0}],\,|z|>1,
\end{array}
\label{jksq}
\end{equation}
and for $|z|=1$,
\begin{equation}
\begin{array}{rl}
\widehat{\Pi}_{1}(z)=&e^{\Gamma_{1}^{+}(z)}[\Psi^{+}(z)+c_{1}z+c_{0}],\\
\widehat{\Pi}_{2}(z)=&z^{-1}e^{\Gamma_{1}^{-}(z)}[\Psi^{-}(z)+c_{1}z+c_{0}],
\end{array}
\label{hjsk}
\end{equation}
where $c_{0}$, $c_{1}$ are constants to be specified,
\begin{displaymath}
\begin{array}{rl}
\Gamma_{1}(z)=&\frac{1}{2\pi i}\int_{|\tau|=1}\log G(\tau)\frac{d\tau}{\tau-z},|z|=1,\\
\Psi(z)=&\frac{1}{2\pi i}\int_{|\tau|=1}g(\tau)e^{-\Gamma_{1}^{+}(\tau)}\frac{d\tau}{\tau-z},|z|=1,
\end{array}
\end{displaymath}
and for $|z_{0}|=1,$
\begin{displaymath}
\begin{array}{rl}
\Gamma_{1}^{+}(z_{0})=\lim_{|z|<1,z\rightarrow z_{0}}\Gamma_{1}(z),&\Psi^{+}(z_{0})=\lim_{|z|<1,z\rightarrow z_{0}}\Psi(z),\\
\Gamma_{1}^{-}(z_{0})=\lim_{|z|>1,z\rightarrow z_{0}}\Gamma_{1}(z),&\Psi^{-}(z_{0})=\lim_{|z|>1,z\rightarrow z_{0}}\Psi(z).
\end{array}
\end{displaymath}
The constants $c_{0}$, $c_{1}$ are obtained from condition 3 above, for $z=0$, $|z|\rightarrow \infty$ by the system
\begin{equation}e^{\Gamma_{1}(0)}[\Psi(0)+c_{0}]=1,\,c_{1}=1.
\label{df}
\end{equation}

Thus, it remains to determine $\Pi(0,0)$. Combining $\Pi(1,0)/\Pi(0,0)$, which is determined by (\ref{hjsk}), (\ref{df}), with equation (\ref{df1}) we can determine $\Pi(0,0)$.

Then, $\Pi(f_{1}(z),0)$, $z\in L\cup L^{+}$ and $\Pi(0,f_{2}(z))$, $z\in L\cup L^{-}$ are known. Clearly, $f_{i}(z)$, $i=1,2,$ maps $L^{+}$ conformally onto $L_{i}^{+}$. Then, denote by $z=w_{i}(z_{i})=f^{-1}_{i}(z_{i})$, $z_{i}\in L_{i}^{+}$ the inverse mapping. Thus,
\begin{displaymath}
\begin{array}{rl}
\Pi_{1}(z_{1},0)=&\Pi(0,0)e^{\Gamma_{1}(w_{1}(z_{1}))}[\Psi(w_{1}(z_{1}))+c_{1}w_{1}(z_{1})+c_{0}],\,z_{1}\in L_{1}^{+},\vspace{2mm}\\
\Pi_{2}(0,z_{2})=&\Pi(0,0)w^{-1}_{2}(z_{2})e^{\Gamma_{1}(w_{2}(z_{2}))}[\Psi(w_{2}(z_{2}))+c_{1}w_{2}(z_{2})+c_{0}],\,z_{2}\in L_{2}^{+},
\end{array}
\end{displaymath}
Consequently, $\Pi(z_{1},z_{2})$ is determined by the functional equation for $z_{1}\in L_{1}\cup L_{1}^{+}$, $z_{2}\in L_{2}\cup L_{2}^{+}$.
\subsection{Reduction to a Fredholm integral equation of the second kind}\label{fred}
Another approach to cope with the solution of (\ref{equ}) is to reducing it to
a Fredholm integral equation. Indeed, for the contours $L_{1}$, $L_{2}$, defined in (\ref{ww1}) there exists a one-to-one map such as $z_{1}=\omega_{1}(z_{2}):L_{2}\to L_{1}$, $z_{2}=\omega_{2}(z_{1}):L_{1}\to L_{2}$. For $z_{1}\in L_{1}$, $(z_{1},\omega_{2}(z_{1}))$ is a zero pair of the kernel. Let
\begin{displaymath}
\begin{array}{c}
\Omega_{1}(z_{1})=\frac{\Pi(z_{1},0)-\Pi(0,0)}{\Pi(0,0)},\,\Omega_{2}(z_{2})=\frac{\Pi(0,z_{2})-\Pi(0,0)}{z_{2}\Pi(0,0)}.
\end{array}
\end{displaymath}
Note that $\Omega_{1}(0)=0$, $\Omega_{2}(0)=\Pi^{-1}(0,0)\frac{d}{dz}\Pi(0,z_{2})|_{z_{2}=0}$. For $z_{1}\in L_{1}$, $z_{2}=\omega_{2}(z_{1})$, let also
\begin{displaymath}
\begin{array}{c}
T(z_{1})=\frac{\tilde{A}(z_{1},\omega_{2}(z_{1}))-z_{1}}{z_{1}[\omega_{2}(z_{1})-\tilde{B}(z_{1},\omega_{2}(z_{1}))]},\,\,t(z_{1})=\frac{\beta_{3}^{*}(y)-1}{\omega_{2}(z_{1})-\tilde{B}(z_{1},\omega_{2}(z_{1}))}.
\end{array}
\end{displaymath}
Using results from subsection \ref{rh}, $ind_{z_{1}\in L_{1}}T(z_{1})=0$. The functional equation (\ref{equ}) is now rewritten as,
\begin{equation}
\begin{array}{c}
\Omega_{2}(\omega_{2}(z_{1}))=T(z_{1})\Omega_{1}(z_{1})+t(z_{1}),\,z_{1}\in L_{1}.
\end{array}
\label{bn}
\end{equation}
Since $\Omega_{1}(z_{1})$ is regular for $z_{1}\in L_{1}^{+}$ and continuous for $z_{1}\in L_{1}^{+}\cup L_{1}$, and similarly, $\Omega_{2}(z_{2})$ is regular for $z_{2}\in L_{2}^{+}$ and continuous for $z_{2}\in L_{2}^{+}\cup L_{2}$, we have that
\begin{eqnarray}
\begin{array}{c}
\frac{1}{2}\Omega_{k}(z_{k})=\frac{1}{2\pi i}\int_{z_{k}\in L_{k}}\frac{\Omega_{k}(z)}{z-z_{k}}dz,\,z_{k}\in L_{k},\,k=1,2.
\end{array}
\label{ml}
\end{eqnarray}
Substituting (\ref{bn}) in (\ref{ml}) for $k=2$, we arrive after some algebra in,
\begin{displaymath}
\begin{array}{rl}
\frac{1}{2}\Omega_{1}(\omega_{1}(z_{2}))T(\omega_{1}(z_{2}))=&\frac{1}{2\pi i}\int_{z_{2}\in L_{2}}\Omega_{1}(\omega_{1}(z))T(\omega_{1}(z))\frac{dz}{z-z_{2}}\\
&+\frac{1}{2\pi i}\int_{z_{2}\in L_{2}}[t(\omega_{1}(z))-t(\omega_{1}(z_{2}))]\frac{dz}{z-z_{2}},\,z_{2}\in L_{2}.
\end{array}
\end{displaymath}
By substituting $\omega_{1}(z_{2})=z_{1}$, and noticing that when $z_{1}$ traverses $L_{1}$ counterclockwise, then, $z_{2}=\omega_{2}(z_{1})$ traverses $L_{2}$ clockwise, it follows for $z_{1}\in L_{1}$,
\begin{equation}
\begin{array}{rl}
\frac{1}{2}\Omega_{1}(z_{1})=&-\frac{1}{2\pi i}\int_{z_{1}\in L_{1}}\frac{\Omega_{1}(z)T(z)}{T(z_{1})}\frac{\omega_{2}^{\prime}(z)}{\omega_{2}^(z)-\omega_{2}(z_{1})}dz\\
&-\frac{1}{T(z_{1})2\pi i}\int_{z_{1}\in L_{1}}[t(z)-t(z_{1})]\frac{\omega_{2}^{\prime}(z)}{\omega_{2}^(z)-\omega_{2}(z_{1})}dz,
\end{array}
\label{rww}
\end{equation}
for $T(z_{1})\neq 0$, $z_{1}\in L_{1}$. Using (\ref{rww}), the fact that $\Omega_{1}(0)= 0$,  the regularity of $\Omega_{1}(z_{1})$, $z_{1}\in L_{1}^{+}$, and (\ref{ml}) for $k=1$, we have that for $z_{1}\in L_{1}$,
\begin{equation}
\begin{array}{rl}
\Omega_{1}(z_{1})=&\frac{1}{2\pi i}\int_{z_{1}\in L_{1}}\Omega_{1}(z)[\frac{1}{z-z_{1}}-\frac{T(z)}{T(z_{1})}\frac{\omega_{2}^{\prime}(z)}{\omega_{2}^(z)-\omega_{2}(z_{1})}-\frac{1}{z}]dz\\
&-\frac{1}{T(z_{1})2\pi i}\int_{z_{1}\in L_{1}}[t(z)-t(z_{1})]\frac{\omega_{2}^{\prime}(z)}{\omega_{2}^(z)-\omega_{2}(z_{1})}dz,
\end{array}
\label{rww1}
\end{equation}
which is a non-homogeneous Fredholm integral equation of
the second kind for $\Omega_{1}(z_{1})$, $z_{1}\in L_{1}$ \cite{ga} and since $ind_{z_{1}\in L_{1}}T(z_{1})=0$ it has a unique solution, which is the boundary value of a function regular in $L_{1}^{+}$.
There are standard techniques to solve (\ref{rww1}) numerically. Just to mention that the numerical evaluation using the approach of reducing (\ref{equ}) into a Fredholm integral equation requires
the lesser computational effort.
\subsection{Basic performance metrics and numerical issues}
In the following we derive expressions for the mean number of customers in orbits at a departure instant. Substituting $z_{2}=1$, and using the  in the functional equation yields
\begin{displaymath}
\begin{array}{rl}
E(X_{1})=&\frac{d}{dz_{1}}\Pi(z_{1},1)|_{z_{1}=1}=W_{1}\Pi(1,0)+W_{2}\Pi(0,1)+W_{0}\Pi(0,0)\\
&+\frac{\theta(\lambda^{2}(\overline{b}_{3}-\overline{b})-\theta)}{2(\lambda+\theta)(2\lambda+\theta)}\frac{d}{dz_{1}}\Pi(z_{1},0)|_{z_{1}=1},
\end{array}
\end{displaymath}
where,
\begin{displaymath}
\begin{array}{rl}
W_{0}=&\frac{\theta\lambda(2\overline{b}_{3}-\overline{b}+\frac{\lambda}{2}(\overline{b}_{3}^{(2)}-\overline{b}^{(2)}))}{(2\lambda+\theta)(1-\rho)}
-\frac{KH(\lambda+\theta)}{(1-\rho)^{2}}\frac{1+\lambda(\overline{b}_{3}-\overline{b})}{2\lambda+\theta},\\
W_{1}=&\frac{\lambda(2\lambda(\overline{b}_{3}-\overline{b})-\theta\overline{b}+\frac{\lambda^{2}}{2}(\overline{b}_{3}^{(2)}-\overline{b}^{(2)}))}{(2\lambda+\theta)(1-\rho)}
-\frac{KH(\lambda+\theta)}{\theta(1-\rho)^{2}}\frac{\lambda^{2}(\overline{b}_{3}-\overline{b})-\theta}{2\lambda+\theta},\\
W_{2}=&\frac{\lambda(2\lambda\overline{b}_{3}+\theta\overline{b}+\frac{\lambda^{2}}{2}(\overline{b}_{3}^{(2)}-\overline{b}^{(2)}))}{(2\lambda+\theta)(1-\rho)}
-\frac{KH(\lambda+\theta)}{(1-\rho)^{2}}\frac{\theta+\lambda(2+\lambda(\overline{b}_{3}-\overline{b}))}{2\lambda+\theta},\\
KH=&-\frac{\lambda[\frac{\lambda}{2}(\theta\overline{b}^{(2)}+\lambda\overline{b}_{3}^{(2)})+\theta\overline{b}+2\lambda\overline{b}_{3}]}{2(\lambda+\theta)},
\end{array}
\end{displaymath}
and
\begin{displaymath}
\begin{array}{l}
\frac{d}{dz_{1}}\Pi(z_{1},0)|_{z_{1}=1}=\Pi(0,0)\vspace{2mm}\\
\times \left\{\begin{array}{ll} \frac{d}{dz_{1}}\left(e^{\Gamma_{1}^{+}(w_{1}(z_{1}))}[\Psi^{+}(w_{1}(z_{1}))+c_{1}w_{1}(z_{1})+c_{0}]\right)|_{z_{1}=1}, & if\ z_{1}=1\in L_{1},\\
 \frac{d}{dz_{1}}\left(AC(e^{\Gamma_{1}(w_{1}(z_{1}))}[\Psi(w_{1}(z_{1}))+c_{1}w_{1}(z_{1})+c_{0}])\right)|_{z_{1}=1}, & if\ z_{1}=1\in L_{1}^{-},
\end{array}\right.
\end{array}
\end{displaymath}
where $AC(F(z))$ represents the analytic continuation of a function $F(z)$. When $z_{1}=1\in L_{1}$,
\begin{displaymath}
\begin{array}{rl}
\frac{d}{dz_{1}}\Pi(z_{1},0)|_{z_{1}=1}=&\Pi(0,0)\lim_{z_{1}\in L_{1}^{+},z_{1} \to \infty}\frac{d}{dz_{1}}[(e^{\Gamma_{1}(w_{1}(z_{1}))}[\Psi(w_{1}(z_{1}))+c_{1}w_{1}(z_{1})+c_{0}]]\\
=&\lim_{z_{1}\in L_{1}^{+},z_{1} \to \infty}\left\{w_{1}^{(1)}(z_{1})e^{\Gamma_{1}(z_{1})}\left([\Psi(w_{1}(z_{1}))+w_{1}(z_{1})+c_{0}]\right.\right.\\
&\left.\left.\times\frac{1}{2\pi i}\int_{\tau\in L_{1}}\frac{\log G(\tau)d\tau}{(\tau-w_{1}(z_{1}))^{2}}+\frac{1}{2\pi i}\int_{\tau\in L_{1}}e^{-\Gamma_{1}^{(+)}(\tau)}\frac{\log g(\tau)d\tau}{(\tau-w_{1}(z_{1}))^{2}}\right)+1\right\},
 \end{array}
\end{displaymath}
where $w_{k}^{(1)}(z_{1})=\frac{dw_{k}(z_{k})}{dz_{k}}$, $k=1,2$. Analogous calculations can be made for $\frac{d}{dz_{2}}\Pi(0,z_{2})|_{z_{2}=1}$, which lead to the derivation of $E(X_{2})$.

The solution of (\ref{equ}) as described in subsection \ref{rh} is based on
the properties of th conformal mappings $f_{1}:L^{+}\to L_{1}^{+}$, $f_{2}:L^{-}\to L_{2}^{-}$. In equations (\ref{ww2}), (\ref{ww3}) we derived  integral
expressions for these mappings. These expressions contain the
function $\psi(z)$, $z\in L$, which is determined as the solution of an integral equation (\ref{lp}). Such an integral equation cannot be solved explicitly, but numerically. 

There are a lot of existing techniques to solve numerically such integrals (e.g. trapezoid rule), and standard iteration procedures show
rapid convergence based on the values of the parameters. For a detailed treatment of how you can treat numerically (\ref{ww2}), (\ref{ww3}), (\ref{lp}), see \cite{bv}, Ch. IV.2. 

Clearly, the numerical computation of the exact conformal mappings is generally time consuming. Since $L_{1}$, $L_{2}$ are close to ellipses, alternatively, we can approximate them by conformal mappings that map the interior of ellipses to $L^{+}$ \cite{nehar}. In particular, we can approximate the contour $L_{j}$ by ellipse $E_{j}$ with semi-axes $\rho(0)$, $\rho(\pi/2)$, $j=1,2$. Then, $\epsilon(z_{j})$ maps $E_{j}^{+}$ to $L^{+}$ \cite{nehar}, where
\begin{displaymath}
\begin{array}{rl}
\epsilon(z_{j})=\sqrt{k}sn\left(\frac{2Q}{\pi}\sin^{-1}(\frac{z_{j}}{\sqrt{\rho^{2}(0)-\rho^{2}(\pi/2)}});k^{2}\right),&k=16q\prod_{n=1}^{\infty}\left(\frac{1+q^{2n}}{1+q^{2n-1}}\right)^{8},\\
q=\left(\frac{\rho(0)-\rho(\pi/2)}{\rho(0)+\rho(\pi/2)}\right)^{2},&Q=\int_{0}^{1}\frac{dt}{\sqrt{(1+t^{2})(1-k^{2}t^{2})}},
\end{array}
\end{displaymath} 
where $sn(w;l)$ is the Jacobian elliptic function. Our approximation for $f^{-1}_{j}(z_{j})$ is $\epsilon(z_{j})$, $z_{j}\in L_{j}\cup L_{j}^{+}$.
\section{The asymmetrical system}\label{asym}
In the following we investigate the general asymmetric model where $B_{1}\nsim B_{2}\nsim B_{3}$, $\lambda_{1}\neq\lambda_{2}$, $\theta_{1}\neq\theta_{2}$. We proceed with the analysis of the kernel.
\subsection{Analysis of the kernel}
Since $\Pi(z_{1},z_{2})$ is a generating function, it should be regular for
$\left|z_{1}\right| < 1$, continuous for $\left|z_{1}\right| \leq1$, for every fixed $z_{2}$ with $\left|z_{2}\right| \leq 1$; and similarly,
with $z_{1}$ and $z_{2}$ interchanged. This implies that every zerotuple of the kernel $K(z_{1}, z_{2})$
should be a zerotuple of the right-hand side of (\ref{equ}). Hence, we first need to analyze
the zeros of the kernel
\begin{equation}
\begin{array}{c}
K(z_{1},z_{2})=z_{1}z_{2}-\frac{\theta_{1}}{\lambda+\theta}z_{2}\beta_{1}^{*}(y)-\frac{\theta_{2}}{\lambda+\theta}z_{1}\beta_{2}^{*}(y)-\frac{\lambda z_{1}z_{2}}{\lambda+\theta}\beta_{3}^{*}(y).
\end{array}
\label{kert}
\end{equation}
Let $w_{i}=2r_{i}z_{i},\,i=1,2,\,\delta=\frac{1}{2}(w_{1}+w_{2})$. Then, $K(w_{1}/2r_{1},w_{2}/2r_{2})$ is well defined for $w_{1},\,w_{2}$ with $\Re(\delta)\leq 1$, and
$K(w_{1}/2r_{1},(2\delta-w_{1})/2r_{2})$ is a quadratic equation in $w_{1}$ for every fixed $\delta$ with $\Re(\delta)\leq 1$.
After some algebra, the equation $4r_{1}r_{2}K(w_{1}/2r_{1},(2\delta-w_{1})/2r_{2})=0$ can be written as:
\begin{equation}
\begin{array}{l}
[w_{1}-\left(\delta+\frac{\theta_{1}r_{1}\beta_{1}^{*}(\lambda(1-\delta))-\theta_{2}r_{2}\beta_{2}^{*}(\lambda(1-\delta))}{\theta+\lambda(1-\beta_{3}^{*}(\lambda(1-\delta)))}\right)]^{2}\\=[\delta-\frac{\theta_{1}r_{1}\beta_{1}^{*}(\lambda(1-\delta))+\theta_{2}r_{2}\beta_{2}^{*}(\lambda(1-\delta))}{\theta+\lambda(1-\beta_{3}^{*}(\lambda(1-\delta)))}]^{2}-\frac{4\theta_{1}r_{1}\theta_{2}r_{2}\beta_{1}^{*}(\lambda(1-\delta))\beta_{2}^{*}(\lambda(1-\delta))}{(\theta+\lambda(1-\beta_{3}^{*}(\lambda(1-\delta))))^{2}}.
\end{array}
\label{del}
\end{equation}
It must be noted that the right hand side of (\ref{del}) is the discriminant of the quadratic equation $4r_{1}r_{2}K(w_{1}/2r_{1},(2\delta-w_{1})/2r_{2})=0$. Introduce the following function for $\Re(\delta)\leq1,\,0\leq\phi\leq2\pi$
\begin{equation}
\begin{array}{c}
k(\delta)=\delta-\frac{\theta_{1}r_{1}\beta_{1}^{*}(\lambda(1-\delta))+\theta_{2}r_{2}\beta_{2}^{*}(\lambda(1-\delta))}{\theta+\lambda(1-\beta_{3}^{*}(\lambda(1-\delta)))}-\frac{2\cos\phi\sqrt{\theta_{1}r_{1}\theta_{2}r_{2}}\sqrt{\beta_{1}^{*}(\lambda(1-\delta))\beta_{2}^{*}(\lambda(1-\delta))}}{\theta+\lambda(1-\beta_{3}^{*}(\lambda(1-\delta)))}.
\end{array}
\label{po}
\end{equation}
\begin{lemma}\label{lem1}
For every $\phi\in[0,2\pi]$, the function $k(\delta)$ has exactly one real zero on every
one of its two branches, say $\delta_{i}=\delta_{i}(\phi)$, of multiplicity one in $\Re(\delta)\leq 1$. Furthermore,
$\delta_{1}(\phi)=\delta_{2}(\phi+\pi)$.
\end{lemma}
\begin{proof}
See Appendix.
\end{proof}
Substitution of $\delta(\phi)\equiv\delta_{1}(\phi)$ into (\ref{del}) yields
\begin{equation}
\begin{array}{rl}
w_{1}(\phi)=&\frac{2}{\Delta(\phi)}\left[\theta_{1}r_{1}\beta_{1}^{*}(\lambda(1-\delta(\phi)))\right.\\&\left.+e^{i\phi}\sqrt{r_{1}\theta_{1}r_{2}\theta_{2}}\sqrt{\beta_{1}^{*}(\lambda(1-\delta(\phi)))\beta_{2}^{*}(\lambda(1-\delta(\phi)))}\right],\vspace{2mm}\\
w_{2}(\phi)=&\frac{2}{\Delta(\phi)}\left[\theta_{2}r_{2}\beta_{2}^{*}(\lambda(1-\delta(\phi)))\right.\\&\left.+e^{-i\phi}\sqrt{r_{1}\theta_{1}r_{2}\theta_{2}}\sqrt{\beta_{1}^{*}(\lambda(1-\delta(\phi)))\beta_{2}^{*}(\lambda(1-\delta(\phi)))}\right],
\end{array}
\label{ww}
\end{equation}
where $\Delta(\phi)=\theta+\lambda(1-\beta_{3}^{*}(\lambda(1-\delta(\phi)))).$

Equation (\ref{ww}) defines a one-to-one mapping $f$ between $w_{1}(\phi)$ and $w_{2}(\phi)$, i.e. $w_{1}(\phi)=f(w_{2}(\phi))$ or $w_{2}(\phi)=f^{-1}(w_{1}(\phi))$, where $f^{-1}$ denotes the inverse of $f$. For $\phi\in[0,2\pi]$, put
\begin{displaymath}
\begin{array}{rl}
\rho_{1}(\phi)=\left|w_{1}(\phi)\right|,\,u_{1}(\phi)=\arg w_{1}(\phi),&\rho_{2}(\phi)=\left|w_{2}(\phi)\right|,\,u_{2}(\phi)=-\arg w_{2}(\phi).
\end{array}
\end{displaymath}
Then,
\begin{displaymath}
\begin{array}{c}
w_{1}(\phi)=\rho_{1}(\phi)e^{i u_{1}(\phi)},\,\,\,w_{2}(\phi)=\rho_{2}(\phi)e^{-i u_{2}(\phi)},\,\phi\in[0,2\pi].
\end{array}
\end{displaymath}
Let $L_{1}=\left\{w_{1}=w_{1}(\phi);\ \phi\in[0,2\pi]\right\},\,\,\,\,\,\,L_{2}=\left\{w_{2}=w_{2}(\phi);\ \phi\in[0,2\pi]\right\}.$
\begin{lemma}\label{lem2}
$L_{1}$, $L_{2}$ are simple smooth contours.
\end{lemma}
\begin{proof}
See Appendix.
\end{proof}

The contours $L_{1}$, $L_{2}$ satisfy the conditions given in \cite{coh2}. Therefore, using Theorem 1.1 in \cite{coh2} there exists a unique simple contour $L$ in the $w$-plane with $w=0\in L^{+},\,\,\,w=1\in L$, $w=\infty \in L^{-}$, and functions, $f_{1}(w):L^{+}\cup L\rightarrow L_{1}^{+}\cup L_{1},\,\,\,f_{2}(w):L^{-}\cup L\rightarrow L_{2}^{+}\cup L_{2}$, such that: $i)$ $w=0$ is a simple zero of $f_{1}(.)$, and $w=\infty$ is a simple zero of $f_{2}(.)$,
$ii)$ $f_{1}:L^{+}\rightarrow L_{1}^{+}$ is regular and univalent for $w\in L^{+}$, $iii)$ $f_{2}:L^{-}\rightarrow L_{2}^{+}$ is regular and univalent for $w\in L^{-}$, $iv)$ $f_{1}(w)=f(f_{2}(w))$, $w\in L$, where $L^{+}$ denotes the interior of the contour $L$, and $L^{-}$ its exterior. 

These conformal mappings can be constructed as the solution of a simple Riemann boundary value problem (see \cite{bv}, pp. 92-100),
\begin{displaymath}
\begin{array}{rl}
f_{1}(w)=&w\frac{\rho(\psi(1))}{\rho_{1}(\psi(1))}\exp\left(\frac{1}{2\pi i}\int_{\zeta\in L}[\log\rho(\psi(\zeta))+iu(\psi(\zeta))](\frac{\zeta+w}{\zeta-w}-\frac{\zeta+1}{\zeta-1})\frac{d\zeta}{\zeta}\right),\,w\in L^{+},\\
f_{2}(w)=&w^{-1}\frac{\rho(\psi(1))}{\rho_{2}(\psi(1))}\exp\left(\frac{-1}{2\pi i}\int_{\zeta\in L}[\log\rho(\psi(\zeta))+iu(\psi(\zeta))](\frac{\zeta+w}{\zeta-w}-\frac{\zeta+1}{\zeta-1})\frac{d\zeta}{\zeta}\right),w\in L^{-},
\end{array}
\end{displaymath}
and for $w\in L$,
\begin{equation}
\begin{array}{rl}
f_{1}(w)=&w\frac{\rho(\psi(1))}{\rho_{1}(\psi(1))}\rho(\psi(w))\exp\left(\frac{1}{2\pi i}\int_{\zeta\in L}[\log\rho(\psi(\zeta))+iu(\psi(\zeta))](\frac{\zeta+w}{\zeta-w}-\frac{\zeta+1}{\zeta-1})\frac{d\zeta}{\zeta}\right),\\
f_{2}(w)=&w^{-1}\frac{\rho(\psi(1))}{\rho_{2}(\psi(1))}\rho(\psi(w))\exp\left(\frac{-1}{2\pi i}\int_{\zeta\in L}[\log\rho(\psi(\zeta))+iu(\psi(\zeta))](\frac{\zeta+w}{\zeta-w}-\frac{\zeta+1}{\zeta-1})\frac{d\zeta}{\zeta}\right),
\end{array}
\label{qe}
\end{equation}
where,
\begin{displaymath}
\begin{array}{c}
\rho(\psi(w))=\sqrt{\rho_{1}(\psi(w))\rho_{2}(\psi(w))},\,u(\psi(w))=\frac{u_{1}(\psi(w))-u_{2}(\psi(w))}{2},
\end{array}
\end{displaymath}
and $\psi(w)$, $w\in L$, $\psi(1)=0$ is a real function such that $\phi=\psi(w)$, $w\in L$ and
\begin{displaymath}
\begin{array}{rl}
\rho_{1}(\psi(1))=\left|w_{1}(0)\right|=&\frac{2}{\Delta(0)}\left[\theta_{1}r_{1}\beta_{1}^{*}(\lambda(1-\delta(0)))\right.\\&\left.+\sqrt{r_{1}\theta_{1}r_{2}\theta_{2}}\sqrt{\beta_{1}^{*}(\lambda(1-\delta(0)))\beta_{2}^{*}(\lambda(1-\delta(0)))}\right],\vspace{2mm}\\
\rho_{2}(\psi(1))=\left|w_{2}(0)\right|=&\frac{2}{\Delta(0)}\left[\theta_{2}r_{2}\beta_{2}^{*}(\lambda(1-\delta(0)))\right.\\&\left.+\sqrt{r_{1}\theta_{1}r_{2}\theta_{2}}\sqrt{\beta_{1}^{*}(\lambda(1-\delta(0)))\beta_{2}^{*}(\lambda(1-\delta(0)))}\right].
\end{array}
\end{displaymath}
Using (\ref{qe}) the function $\phi=\psi(w)$ is uniquely determined by the following equation for $w\in L$:
\begin{displaymath}
\begin{array}{l}
\left(\frac{\rho_{1}(\psi(w))}{\rho_{2}(\psi(w))}\right)^{1/2}\exp\left(i\frac{u_{1}(\psi(w))+u_{2}(\psi(w))}{2}\right)\\=w\frac{\rho(\psi(1))}{\rho_{1}(\psi(1))}\exp\left(\frac{1}{2\pi i}\int_{\zeta\in L}[\log\rho(\psi(\zeta))+iu(\psi(\zeta))](\frac{\zeta+w}{\zeta-w}-\frac{\zeta+1}{\zeta-1})\frac{d\zeta}{\zeta}\right).
\end{array}
\end{displaymath}

Let now $z_{1}(w)=\frac{f_{1}(w)}{2r_{1}}$, $w\in L^{+}\cup L$, $z_{2}(w)=\frac{f_{2}(w)}{2r_{2}}$ for $w\in L^{-}\cup L$ and
\begin{displaymath}
C_{1}=\left\{z_{1}=z_{1}(w),w\in L\right\},\,C_{2}=\left\{z_{2}=z_{2}(w),w\in L\right\}.
\end{displaymath}

Clearly the pairs $(z_{1}(w),z_{2}(w))$, $w\in L$ are zeros of the kernel $K(z_{1},z_{2})$. Before formulating a Riemann boundary value problem for the functional equation (\ref{equ}), and deriving its solution by using the zerotuples $(z_{1}(w),z_{2}(w))$, $w\in L$, of the kernel $K(z_{1},z_{2})$, we need to carefully check their positions, because with the different choice of parameters, $z_{i}(w)$ for $w\in L$ would be inside, on or outside the unit circle. In the latter case, the analytic continuation for the function of the right-hand side in (\ref{equ}) is necessary. From the construction of $z_{i}(·)$, $i = 1,2$, it can be seen that $z_{i}(w)$ has its maximum modulo at $w = 1$. Since $\psi(1) = 0$, i.e., $w = 1$ corresponds to $\phi = 0$, we have
\begin{displaymath}
\begin{array}{rl}
z_{1}(1)=&\frac{1}{\Delta(0)}[\theta_{1}\beta^{*}_{1}(\lambda(1-\delta(0)))+\frac{\sqrt{\theta_{1}r_{1}\theta_{2}r_{2}}}{r_{1}}\sqrt{\widehat{\beta}^{*}(\lambda(1-\delta(0)))}],\\
z_{2}(1)=&\frac{1}{\Delta(0)}[\theta_{2}\beta^{*}_{2}(\lambda(1-\delta(0)))+\frac{\sqrt{\theta_{1}r_{1}\theta_{2}r_{2}}}{r_{2}}\sqrt{\widehat{\beta}^{*}(\lambda(1-\delta(0)))}],
\end{array}
\end{displaymath}
where $\widehat{\beta}^{*}(\lambda(1-\delta(0)))=\beta_{1}^{*}(\lambda(1-\delta(0)))\beta_{2}^{*}(\lambda(1-\delta(0)))$.

Without loss of generality we assume that $\theta_{2}r_{1}\geq\theta_{1}r_{2}$. Then,
\begin{enumerate}
\item $\theta_{2}r_{1}=\theta_{1}r_{2}$. In such a case, since $\delta(0)=1$,
\begin{displaymath}
z_{i}(1)=\frac{1}{\theta}[\theta_{i}+\frac{\sqrt{\theta_{1}r_{1}\theta_{2}r_{2}}}{r_{i}}]=\frac{\theta_{1}+\theta_{2}}{\theta}=1,\,i=1,2.
\end{displaymath}
\item $\theta_{2}r_{1}>\theta_{1}r_{2}$. Since $\delta(0)<1$, in this case we have
\begin{displaymath}
\begin{array}{rl}
z_{1}(1)=&\frac{1}{\Delta(0)}[\theta_{1}\beta^{*}_{1}(\lambda(1-\delta(0)))+\frac{\sqrt{\theta_{1}r_{1}\theta_{2}r_{2}}}{r_{1}}\sqrt{\widehat{\beta}^{*}(\lambda(1-\delta(0)))}]\\
<&\frac{1}{\theta}[\theta_{1}+\frac{\sqrt{\theta_{1}r_{1}\theta_{2}r_{2}}}{r_{1}}]<\frac{1}{\theta}[\frac{\theta_{2}r_{1}+\theta_{1}r_{2}}{2r_{1}}+\theta_{1}]<\frac{\theta_{1}+\theta_{2}}{\theta}=1.
\end{array}
\end{displaymath}

Let $\widehat{\theta}_{i}=\frac{\theta_{i}}{\theta}$ and $c=(1-\widehat{\theta}_{2}\beta^{*}_{2}(\lambda(1-\delta(0))))\Delta(0)/\sqrt{\widehat{\beta}^{*}(\lambda(1-\delta(0)))}$. Then, $z_{2}(1)<1$, if $r_{1}\widehat{\theta}_{2}<r_{2}c^{2}/\widehat{\theta}_{1}$, $z_{2}=1$, if $r_{1}\widehat{\theta}_{2}=r_{2}c^{2}/\widehat{\theta}_{1}$ and $z_{2}(1)>1$ if $r_{1}\widehat{\theta}_{2}>r_{2}c^{2}/\widehat{\theta}_{1}$.
\end{enumerate}
\subsection{Solution of the functional equation}
In this section we formulate a Riemann boundary value problem for the functional equation (\ref{equ}), and derive its solution by using the zerotuple $(z_{1}(w),z_{2}(w))$, $w\in L$, of the kernel $K(z_{1},z_{2})$. Since $\Pi(z_{1},z_{2})$ should be regular for $|z_{1}| < 1$, continuous for $|z_{1}|\leq 1$, for every fixed $z_{2}$ with $|z_{1}|\leq 1$; and similarly, with $z_{1}$ and $z_{2}$ interchanged, the right-hand side of (\ref{equ}) should be zero for all those $w\in L$, for which $(z_{1}(w),z_{2}(w))$ forms a pair of zeros of $K(z_{1},z_{2})$ inside the product of unit circles. $|z_{1}(w)|\leq 1$ always holds for $w\in L$, but $|z_{2}(w)|\leq 1$ may not hold for $w\in L$. By analytic continuation, we can prove that the right-hand side of (\ref{equ}) also is zero in the case that $z_{2}(w)$ is not inside the unit circle. Hence we have the following relation:
\begin{equation}
\begin{array}{l}
\widehat{\Pi}(z_{1}(w),0)z_{2}(w)[z_{1}(w)-\widetilde{A}(w)]=z_{1}(w)[\widetilde{B}(w)-z_{2}(w)]\widehat{\Pi}(0,z_{2}(w))\\
+[z_{1}(w)z_{2}(w)+z_{2}(w)(r_{1}z_{1}(w)\beta^{*}_{3}(y(w))-\widetilde{A}(w))+z_{1}(w)(r_{2}z_{2}(w)\beta^{*}_{3}(y(w))-\widetilde{B}(w))]
\end{array}
\label{wq}
\end{equation}
where for $w\in L$, $y(w)=\lambda(1-r_{1}z_{1}(w)-r_{2}z_{2}(w))$ and,
\begin{displaymath}
\begin{array}{rl}
\widetilde{A}(w)=\widetilde{A}(z_{1}(w)),z_{2}(w),&\widetilde{B}(w)=\widetilde{B}(z_{1}(w),z_{2}(w)),\vspace{2mm}\\
\widehat{\Pi}(z_{1}(w),0)=\widehat{\Pi}_{1}(w)=\frac{\Pi(z_{1}(w),0)}{\Pi(0,0)},&\widehat{\Pi}(0,z_{2}(w))=\widehat{\Pi}_{2}(w)=\frac{\Pi(0,z_{2}(w))}{\Pi(0,0)}.
\end{array}
\end{displaymath}
Thus, (\ref{wq}) can be written as
\begin{equation}
\widehat{\Pi}_{1}(w)=G(w)\widehat{\Pi}_{2}(w)+g(w),\,w\in L,
\label{wq1}
\end{equation}
where,
\begin{displaymath}
\begin{array}{rl}
G(w)=&-\frac{1-\frac{\widetilde{B}(w)}{z_{2}(w)}}{1-\frac{\widetilde{A}(w)}{z_{1}(w)}}=-\frac{1-(1-a_{2})\beta_{3}^{*}(y(w))}{1-(1-a_{1})\beta_{3}^{*}(y(w))}\frac{1-\frac{\widetilde{\beta}_{2}^{*}(w)}{z_{2}(w)}}{1-\frac{\widetilde{\beta}_{1}^{*}(w)}{z_{1}(w)}},\\
g(w)=&-G(w)+\frac{\beta^{*}_{3}(y(w))-\frac{\widetilde{A}(w)}{z_{1}(w)}}{1-\frac{\widetilde{A}(w)}{z_{1}(w)}}=-G(w)-\frac{1-\beta_{3}^{*}(y(w))}{(1-(1-a_{1})\beta_{3}^{*}(y(w)))[1-\frac{\widetilde{\beta}_{1}^{*}(w)}{z_{1}(w)}]}+1,
\end{array}
\end{displaymath}
where, for $w\in L$, $a_{j}=\frac{\theta_{j}}{\lambda+\theta_{j}}$, $j=1,2$, 
\begin{displaymath}
\begin{array}{c}
\widetilde{\beta}_{j}^{*}(w)=\widetilde{\beta}_{j}^{*}(z_{1}(w),z_{2}(w))=\frac{a_{j}\beta_{j}^{*}(y(w))}{1-(1-a_{j})\beta_{3}^{*}(y(w))}=\frac{\theta_{j}\beta_{j}^{*}(y(w))}{\theta_{j}+\lambda(1-\beta_{3}^{*}(y(w))},
\end{array}
\end{displaymath}

We can easily show that $\widetilde{\beta}_{j}^{*}(z_{1},z_{2})$ has a probabilistic interpretation. Indeed, it is the generating function of the joint orbit queue length distribution of the number of customers that arrive from the epoch a service is initiated until the epoch the server becomes idle for the first time after the service of a customer coming from the orbit queue $j$, $j=1,2$, given that we allow retrials only from orbit queue $j$. Let $S_{j}$ be the corresponding time interval, and denote by $N_{i}(S_{j})$ the number of type $i$ customers that join the orbit queue $i$ during $S_{j}$, $i,j=1,2$. If $\widetilde{b}_{m,n}^{(j)}(t)dt=P(t<S_{j}\leq t+dt,N_{1}(S_{j})=m,N_{2}(S_{j})=n)$, then,
\begin{displaymath}
\begin{array}{c}
\widetilde{b}_{m,n}^{(j)}(t)=a_{j}h_{m,n}(t)b_{j}(t)+(1-a_{j})\sum_{l=0}^{m}\sum_{c=0}^{n}h_{l,c}(t)*\widetilde{b}_{m-l,n-c}^{(j)}(t).
\end{array}
\end{displaymath}
Let ${\beta}_{j}^{*}(z_{1},z_{2},s)=\int_{0}^{\infty}e^{-st}\sum_{m=0}^{\infty}\sum_{n=0}^{\infty}\widetilde{b}_{m,n}^{(j)}(t)z_{1}^{m}z_{2}^{n}dt$. Then, $\widetilde{\beta}_{j}^{*}(z_{1},z_{2})=\widetilde{\beta}_{j}^{*}(z_{1},z_{2},0).$
\begin{theorem}
$\widehat{\Pi}_{1}(w)$ is regular for $w\in L^{+}$, continuous for $w\in L\cup L^{+}$, and $\widehat{\Pi}_{2}(w)$ is regular for $w\in L^{-}$, continuous for $w\in L\cup L^{-}$.
\end{theorem}

Thus, we have the following problem on the unit circle $L$: Find two functions $\widehat{\Pi}_{1}(w)$, $\widehat{\Pi}_{2}(w)$ such that: $i)$ $\widehat{\Pi}_{1}(w)$ is regular for $w\in L^{+}$, continuous for $w\in L\cup L^{+}$, $ii)$ $\widehat{\Pi}_{2}(w)$ is regular for $w\in L^{-}$, continuous for $w\in L\cup L^{-}$, $iii)$ for $w\in L$, $\widehat{\Pi}_{1}(w)=G(w)\widehat{\Pi}_{2}(w)+g(w)$, and
\begin{equation}
\begin{array}{c}
\lim_{|w|\rightarrow \infty}\widehat{\Pi}_{2}(w)=1,\,\widehat{\Pi}_{1}(0)=1,
\end{array}
\label{lpm}
\end{equation}
In order to derive a solution for the Riemann boundary value problem, we need to investigate properties of the functions $G(w)$, $g(w)$, $w\in L$. In particular, $i)$ $0<G(w)<\infty$, $0<g(w)<\infty$, $w\in L$, $ii)$ $G(w)$, $g(w)$ should satisfy the Holder condition on $L$. Instead of them, we will consider the following equivalent functions for $\phi\in[0,2\pi]$:
\begin{displaymath}
\begin{array}{rl}
\widehat{G}(\delta(\phi))=&-\frac{1-(1-a_{2})\beta_{3}^{*}(\lambda(1-\delta(\phi)))}{1-(1-a_{1})\beta_{3}^{*}(\lambda(1-\delta(\phi)))}\frac{1-\frac{\widetilde{\beta}_{2}^{*}(w_{1}(\phi)/2r_{1},w_{2}(\phi)/2r_{2})}{w_{2}(\phi)/2r_{2}}}{1-\frac{\widetilde{\beta}_{1}^{*}(w_{1}(\phi)/2r_{1},w_{2}(\phi)/2r_{2})}{w_{1}(\phi)/2r_{2}}}\vspace{2mm}\\
\widehat{g}(\delta(\phi))=&-\widehat{G}(\delta(\phi))-\frac{1-\beta_{3}^{*}(\lambda(1-\delta(\phi)))}{(1-(1-a_{1})\beta_{3}^{*}(\lambda(1-\delta(\phi))))[1-\frac{\widetilde{\beta}_{1}^{*}(w_{1}(\phi)/2r_{1},w_{2}(\phi)/2r_{2})}{w_{1}(\phi)/2r_{1}}]}+1.
\end{array}
\end{displaymath}

Since $\delta(\phi)$ is real and $L_{1}$, $L_{2}$ are simple contours, the two points $\delta(0)$, $\delta(\pi)$ are the only candidates zeros of the numerator and the denominator of $\widehat{G}(\delta(\phi))$, which takes after some algebra the following form,
\begin{displaymath}
\begin{array}{c}
\widehat{G}(\delta(\phi))=-\frac{r_{2}w_{1}(\phi)(\lambda+\theta_{1})\Delta_{2}(\phi)}{r_{1}w_{2}(\phi)(\lambda+\theta_{2})\Delta_{1}(\phi)}\left(\frac{\frac{w_{2}(\phi)}{2r_{2}}-\frac{\theta_{2}\beta_{2}^{*}(\lambda(1-\delta(\phi)))}{\theta_{2}+\lambda(1-\beta_{3}^{*}(\lambda(1-\delta(\phi))))}}{\frac{w_{1}(\phi)}{2r_{1}}-\frac{\theta_{1}\beta_{1}^{*}(\lambda(1-\delta(\phi)))}{\theta_{1}+\lambda(1-\beta_{3}^{*}(\lambda(1-\delta(\phi))))}}\right),
\end{array}
\end{displaymath}
where $\Delta_{j}(\phi)=\theta_{j}+\lambda(1-\beta_{3}(\lambda(1-\delta(\phi))))$, $j=1,2$. We concentrate on the part in the parenthesis and see that at the point $\delta(\pi)$, its numerator and denominator become respectively,
\begin{displaymath}
\begin{array}{c}
-\left(\frac{\theta_{1}\theta_{2}\beta_{2}(\lambda(1-\delta(\pi)))+\Delta_{2}(\pi)\frac{\sqrt{\theta_{1}\theta_{2}r_{1}r_{2}}}{r_{2}}\sqrt{\beta_{1}(\lambda(1-\delta(\pi)))\beta_{2}(\lambda(1-\delta(\pi)))}}{\Delta(\pi)}\right)<0,\vspace{2mm}\\
-\left(\frac{\theta_{1}\theta_{2}\beta_{1}(\lambda(1-\delta(\pi)))+\Delta_{1}(\pi)\frac{\sqrt{\theta_{1}\theta_{2}r_{1}r_{2}}}{r_{1}}\sqrt{\beta_{1}(\lambda(1-\delta(\pi)))\beta_{2}(\lambda(1-\delta(\pi)))}}{\Delta(\pi)}\right)<0,
\end{array}
\end{displaymath}
and thus, $\delta(\pi)$ is neither zero nor pole of $\widehat{G}(\delta(\phi))$. We will consider now the point $\delta(0)$: For $r_{2}\theta_{1}=r_{1}\theta_{2}$. Since $\delta(0)=1$, we can easily verify that the numerator and the denominator of $\widehat{G}(\delta(\phi))$ vanish simultaneously. Thus, $\delta(0)$ is a cancelled point of $\widehat{G}(\delta(\phi))$ and $0<\left|\widehat{G}(\delta(\phi))\right|<\infty$.

Let now, $r_{1}\theta_{2}>r_{2}\theta_{1}$. Since $\delta(0)<1$, $\beta_{j}(\lambda(1-\delta(0)))<1$, $j=1,2.$ Then,
\begin{enumerate}
\item if $r_{1}\widehat{\theta}_{2}\geq r_{2}c^{2}/\widehat{\theta}_{1}$, then $w_{2}(0)/2r_{2}\geq 1$ and as a result $\frac{w_{2}(0)}{2r_{2}}-\widetilde{\beta}_{2}^{*}(w_{1}(0)/2r_{1},w_{2}(0)/2r_{2})>0$. Note that when $2r_{1}>1$ there is a possibility that the denominator of $\widehat{G}(\delta(\phi))$ vanishes. It is easily seen after some algebra that the denominator of $\widehat{G}(\delta(\phi))$ never vanishes if
\begin{displaymath}
\begin{array}{c}
\frac{r_{2}\beta_{2}^{*}(\lambda(1-\delta(0)))}{r_{1}}\neq \frac{\theta_{1}\theta_{2}\beta_{1}^{*}(\lambda(1-\delta(0)))}{\Delta_{1}^{2}(0)}.
\end{array}
\end{displaymath}
\item if $r_{1}\widehat{\theta}_{2}<r_{2}c^{2}/\widehat{\theta}_{1}$, then $w_{2}(0)/2r_{2}\leq 1$ and in this case we cannot exclude the possibility that the the numerator and the denominator of $\widehat{G}(\delta(0))$ vanish simultaneously. Letting the numerator and denominator equal zero respectively, we obtain the following equalities:
\begin{displaymath}
\begin{array}{rl}
r_{2}\theta_{1}\frac{\theta_{2}\beta_{2}^{*}(\lambda(1-\delta(0)))}{\theta_{2}+\lambda(1-\beta_{3}^{*}(\lambda(1-\delta(0))))}=&\Delta_{1}(0)r_{1}\beta_{1}^{*}(\lambda(1-\delta(0))),\\
r_{1}\theta_{2}\frac{\theta_{1}\beta_{1}^{*}(\lambda(1-\delta(0)))}{\theta_{1}+\lambda(1-\beta_{3}^{*}(\lambda(1-\delta(0))))}=&\Delta_{2}(0)r_{2}\beta_{2}^{*}(\lambda(1-\delta(0))).
\end{array}
\end{displaymath}
\end{enumerate}

Let $\chi$ be the index of the function $G(w)$, $w\in L$,
\begin{displaymath}
\begin{array}{rl}
\chi=&ind_{w\in L}G(w)=ind_{\phi\in[0,2\pi]}\widehat{G}(\delta(\phi))\\
=&ind_{\phi\in[0,2\pi]}\frac{w_{1}(\phi)}{2r_{1}}-ind_{\phi\in[0,2\pi]}\frac{w_{2}(\phi)}{2r_{2}}+ind_{\phi\in[0,2\pi]}[\frac{w_{2}(\phi)}{2r_{2}}-\widetilde{\beta}_{2}^{*}(w_{1}(\phi)/2r_{1},w_{2}(\phi)/2r_{2})]\\
&-ind_{\phi\in[0,2\pi]}[\frac{w_{1}(\phi)}{2r_{1}}-\widetilde{\beta}_{1}^{*}(w_{1}(\phi)/2r_{1},w_{2}(\phi)/2r_{2})].
\end{array}
\end{displaymath}

\begin{lemma}\label{lem3}
Let $r_{2}\theta_{1}\leq r_{1}\theta_{2}$ and $r_{1}\widehat{\theta}_{2}\geq r_{2}c^{2}/\widehat{\theta}_{1}$. Under the assumption
\begin{displaymath}
\begin{array}{c}
\theta{2}r_{1}\widetilde{\beta}_{1}^{*}(w_{1}(0)/2r_{1},w_{2}(0)/2r_{2})>\Delta_{1}(0)r_{2}\beta_{2}^{*}(\lambda(1-\delta(0))),
\end{array}
\end{displaymath}
 the index $\chi=1.$
\end{lemma}
\begin{proof}
See Appendix.
\end{proof}

Then using the standard approach \cite{ga}, we derive the following solution of the non-homogeneous Riemann boundary value problem (\ref{wq1}):
\begin{equation}
\begin{array}{rl}
\widehat{\Pi}_{1}(w)=&e^{\Gamma_{1}(w)}[\Psi(w)+c_{1}w+c_{0}],\,w\in L^{+},\\
\widehat{\Pi}_{2}(w)=&w^{-1}e^{\Gamma_{1}(w)}[\Psi(w)+c_{1}w+c_{0}],\,w\in L^{-},
\end{array}
\label{jk}
\end{equation}
\begin{equation}
\begin{array}{rl}
\widehat{\Pi}_{1}(w)=&e^{\Gamma_{1}^{+}(w)}[\Psi^{+}(w)+c_{1}w+c_{0}],\,w\in L,\\
\widehat{\Pi}_{2}(w)=&w^{-1}e^{\Gamma_{1}^{-}(w)}[\Psi^{-}(w)+c_{1}w+c_{0}],\,w\in L,
\end{array}
\label{hjk}
\end{equation}
where $c_{0}$, $c_{1}$ are constants to be specified from (\ref{jk}), (\ref{lpm}) for $w=0$, $|w|\rightarrow \infty$ by the system
\begin{displaymath}
\begin{array}{c}
e^{\Gamma_{1}(0)}[\Psi(0)+c_{0}]=1,\,c_{1}=1.
\end{array}
\end{displaymath}

Since $z_{1}(w) = f_{1}(w)/2r_{1}$ for $w\in L^{+}\cup L$ and $z_{2}(w) = f_{2}(w)/2r_{2}$ for $w\in L\cup L^{-}$, the existence of the inverse mapping of $f_{i}(w)$ implies that the inverse mapping of $z_{i}(w)$ exists. Let $w_{i}(z_{i}) = f^{-1}_{i} (2r_{i}w)$, where $f^{-1}_{i}$ denotes the inverse mapping of $f_{i}(w)$, then, $w_{1}(z_{1}):C_{1}\cup C_{1}^{+}\rightarrow L\cup L^{+}$, $w_{2}(z_{2}):C_{2}\cup C_{2}^{+}\rightarrow L\cup L^{-}$, are respectively the inverse mappings of $z_{1}(w)$, $z_{2}(w)$. Therefore,
\begin{theorem}
\begin{equation}
\begin{array}{rl}
\widehat{\Pi}_{1}(z_{1})=&e^{\Gamma_{1}(w_{1}(z_{1}))}[\Psi(w_{1}(z_{1}))+c_{1}w_{1}(z_{1})+c_{0}],\,z_{1}\in C_{1}^{+},\\
\widehat{\Pi}_{2}(z_{2})=&(w_{2}(z_{1}))^{-1}e^{\Gamma_{1}(w_{2}(z_{2}))}[\Psi(w_{2}(z_{2}))+c_{1}w_{2}(z_{2})+c_{0}],\,z_{2}\in C_{2}^{+},
\end{array}
\label{jkq}
\end{equation}
and,
\begin{displaymath}
\begin{array}{rl}
\widehat{\Pi}_{1}(z_{1})=&e^{\Gamma_{1}^{+}(w_{1}(z_{1}))}[\Psi^{+}(w_{1}(z_{1}))+c_{1}w_{1}(z_{1})+c_{0}],\,z_{1}\in C_{1},\\
\widehat{\Pi}_{2}(z_{2})=&(w_{2}(z_{2}))^{-1}e^{\Gamma_{1}^{-}(w_{2}(z_{2}))}[\Psi^{-}(w_{2}(z_{1}))+c_{1}w_{2}(z_{2})+c_{0}],\,z_{2}\in C_{2}.
\end{array}
\end{displaymath}
\end{theorem}


\section{Explicit expressions for the completely symmetrical model}\label{sym}
In the following we show how to compute basic performance metrics for the completely symmetrical system without the need of solving a boundary value problem. As a symmetrical model, we mean that $\lambda_{1}=\lambda_{2}=\frac{\lambda}{2}$ (i.e., $r_{1}=r_{2}=\frac{1}{2}$), $\theta_{1}=\theta_{2}=\frac{\theta}{2}$, $B_{j}\sim B$, $j=1,2,3$. Symmetry means also that all queue lengths have the same distributions, but what is more important is that the boundary functions are equal. Using (\ref{equ}), and the fact that $\Pi(1,1)=1$, $\Pi(1,0)=\Pi(0,1)$ we can obtain
\begin{displaymath}
\begin{array}{c}
\Pi(1,0)\frac{2\lambda\theta}{2\lambda+\theta}+\Pi(0,0)\frac{\theta^{2}}{2\lambda+\theta}=\theta-2\rho(\lambda+\theta),
\end{array}
\end{displaymath}
where $\rho=\frac{\lambda}{2}\bar{b}$. Note here that $\theta-2\rho(\lambda+\theta)>0$ due to the stability condition. 

Denote by $\Pi_{1}(z_{1},z_{2})$, $\Pi_{2}(z_{1},z_{2})$, the derivatives of $\Pi(z_{1},z_{2})$ with respect to $z_{1}$, $z_{2}$, respectively. Due to the symmetry, $\Pi_{1}(1,1)=\Pi_{2}(1,1)$, $\Pi_{1}(1,0)=\Pi_{2}(0,1)$. Differentiate (\ref{equ}) with respect to $z_{1}$, and set $(z_{1},z_{2})=(1,1)$ to get
\begin{equation}
\begin{array}{c}
\Pi_{1}(1,1)=\frac{2\rho(\lambda+\theta)(1-\rho)+\frac{\lambda^{2}\bar{b}^{(2)}}{4}(\lambda+\theta)-\frac{\theta^{2}}{2\lambda+\theta}\Pi_{1}(1,0)}{\theta-2\rho(\lambda+\theta)}.
\end{array}
\label{xv}
\end{equation}
Now by setting $z_{1}=z_{2}=z$ in (\ref{equ}) we get,
\begin{equation}
\begin{array}{c}
\frac{d}{dz}\Pi(z,z)|_{z=1}=\frac{2\rho(\lambda+\theta)(1-2\rho)+\frac{\lambda^{2}\bar{b}^{(2)}}{2}(\lambda+\theta)+\frac{2\lambda\theta}{2\lambda+\theta}\Pi_{1}(1,0)}{2(\theta-2\rho(\lambda+\theta))}.
\end{array}
\label{dff}
\end{equation}
However, due to symmetry 
\begin{equation}
\begin{array}{c}
\frac{d}{dz}\Pi(z,z)|_{z=1}=2\Pi_{1}(1,1).
\end{array}
\label{dfgh}
\end{equation}
Substituting (\ref{dfgh}) in (\ref{dff}), and eliminating $\Pi_{1}(1,0)$ from (\ref{dff}) using (\ref{xv}) we obtain
\begin{equation}
\begin{array}{c}
\Pi_{1}(1,1)=\frac{4\rho(2\lambda+\theta-2\rho(\lambda+\theta))+\lambda^{2}\bar{b}^{(2)}(\lambda+\theta)}{4(\theta-2\rho(\lambda+\theta))}.
\end{array}
\label{mn}
\end{equation} 
Since $\Pi_{1}(1,1)$ is equal to the expected number of customers in an orbit, a simple application of Little's law gives the expected orbit delay $E(D)$,
\begin{equation}
\begin{array}{c}
E(D)=\frac{2\Pi_{1}(1,1)}{\lambda}=\frac{4\rho(2\lambda+\theta-2\rho(\lambda+\theta))+\lambda^{2}\bar{b}^{(2)}(\lambda+\theta)}{2\lambda(\theta-2\rho(\lambda+\theta))}.
\end{array}
\label{dela}
\end{equation}
\section{A numerical example}\label{num}
In this section we provide a numerical example regarding the performance of the completely symmetrical system. Assume that the service time is an $Erlang(2,\mu)$ distributed random variable with $\bar{b}=\frac{2}{\mu}$, $\bar{b}^{(2)}=\frac{6}{\mu^{2}}$. 

\begin{figure}[ht!]
\centering
\begin{minipage}[b]{0.45\linewidth}
\includegraphics[scale=0.3]{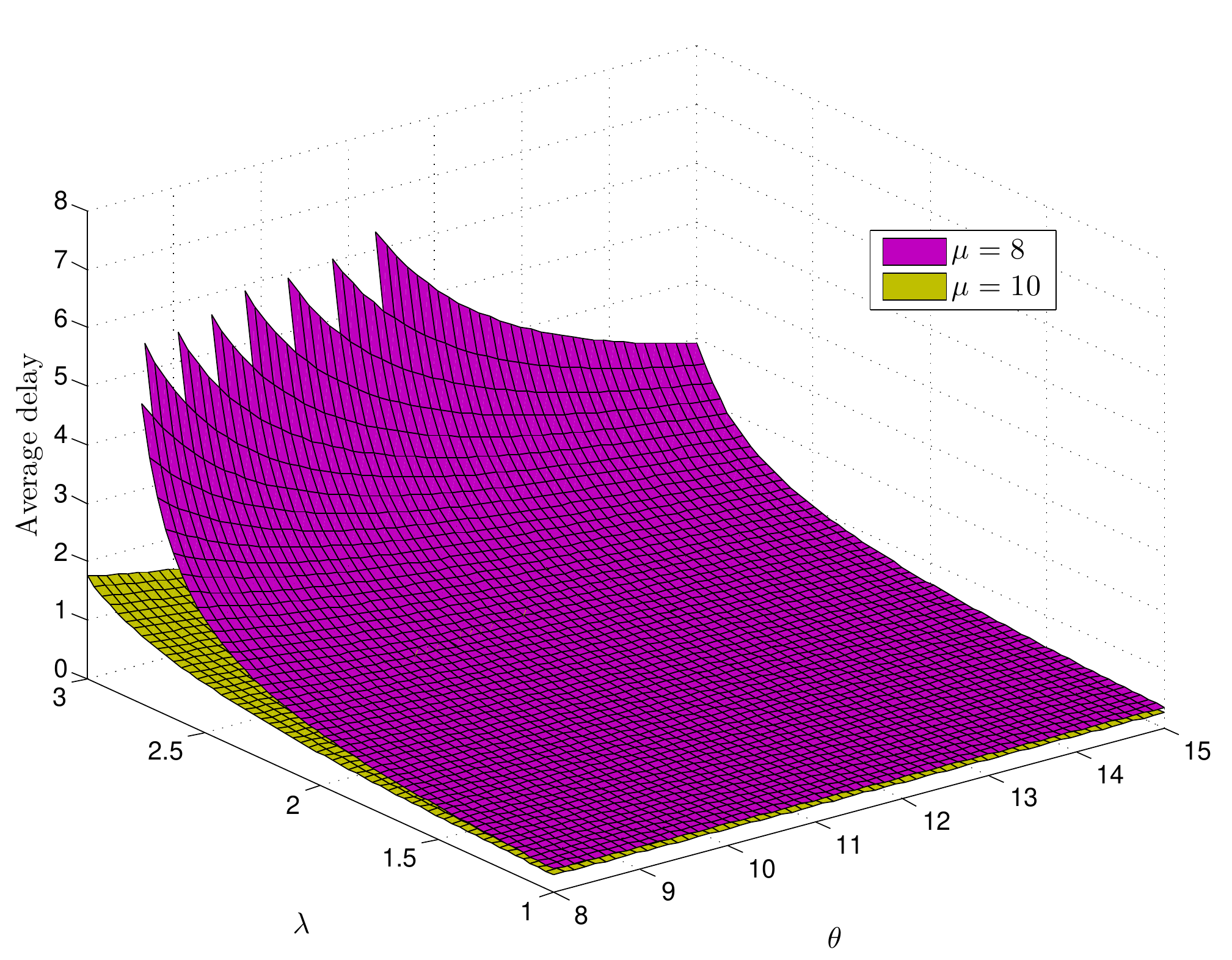}
\end{minipage}
\quad
\begin{minipage}[b]{0.45\linewidth}
\includegraphics[scale=0.3]{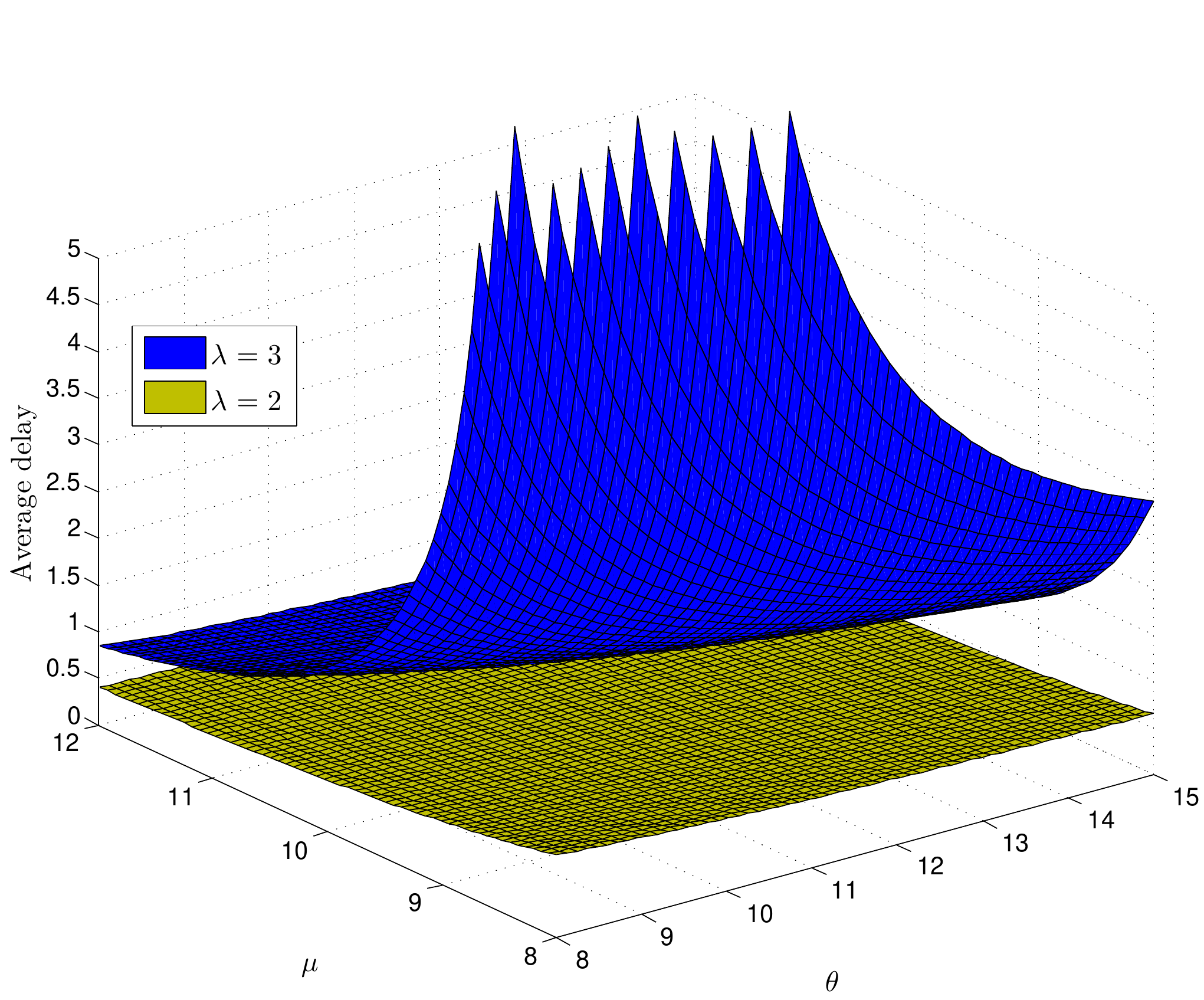}
\end{minipage}\caption{Average delay as a function of $\lambda$, $\theta$ (left), and as a function of $\mu$, $\theta$ (right).}
\label{f1}
\end{figure}

In Figure \ref{f1} (left) we can observe the effect of $\lambda$ and $\theta$ on the average delay obtained in (\ref{dela}). As expected, the increase in $\lambda$ will cause the increase in $E(D)$. That increase becomes more apparent for small values of $\theta$ since in such a case the orbiting customers retry in a ``slow" fashion. Moreover, if the service rate increases, the average delay in an orbit will decrease.

Figure \ref{f1} (right) shows the way $E(D)$ is affected for increasing values of $\mu$, $\theta$. Clearly, the increase in $\mu$ will result in the decrease of the average delay. However, we can easily observe how sensitive is $E(D)$ when we slightly increase $\lambda$, and especially when $\mu$, and $\theta$ take small values. 
\section{Conclusion and future work}\label{conc}
As already mentioned this paper aims to provide a general framework for the fundamental problem of the analysis of multiclass retrial queueing systems with constant retrial policy and general class dependent service times. For the two-orbit case, we generalize the model in \cite{avr}, and provided a compact methodological approach in order to obtain the generating function of the joint orbit queue length distribution in terms of a solution of a Riemann boundary value problem. 

Our results serve as a building block to obtain expressions for the delay in the case of $N$ orbit queues. Clearly, the delay analysis for the general case of $N$ orbit queues under constant retrial policy is highly non-trivial and still remains an open problem. We are currently working towards this direction, and we intent to provide bounds for the queueing delay in an orbit in a
general topology with $N$ orbit queues. Another point of interest is to explore
the possibility to study the heavy traffic behavior of such a model, when the arrivals $\lambda_{j}$ are such that $\widehat{\rho}_{j}\to1$, $j=1,2$. Our approach is also valid for the modelling of even general systems that include vacations, server failures, feedback and arbitrarily distributed retrial times.
\section*{Appendix}
\paragraph{Proof of Lemma \ref{lem1}}
By restricting the function $k(\delta)$ defined in (\ref{po}) to one of its two branches, say
its principal value, we get that $k(\delta)$ is analytic. Let $k(\delta)=\delta-k^{*}(\delta)$, where for $\Re(\delta)\leq1$, $0\leq\phi\leq2\pi$,
\begin{displaymath}
\begin{array}{rl}
k^{*}(\delta)=&\frac{\theta_{1}r_{1}\beta_{1}^{*}(\lambda(1-\delta))+\theta_{2}r_{2}\beta_{2}^{*}(\lambda(1-\delta))}{\theta+\lambda(1-\beta_{3}^{*}(\lambda(1-\delta)))}+\frac{2\cos\phi\sqrt{\theta_{1}r_{1}\theta_{2}r_{2}}\sqrt{\beta_{1}^{*}(\lambda(1-\delta))\beta_{2}^{*}(\lambda(1-\delta))}}{\theta+\lambda(1-\beta_{3}^{*}(\lambda(1-\delta)))}.
\end{array}
\end{displaymath}

If $\theta_{1}r_{2}\neq \theta_{2}r_{1}$, then $2\sqrt{\theta_{1}r_{1}\theta_{2}r_{2}}<\theta_{2}r_{1}+
\theta_{1}r_{2}$. For $\phi\in[0,2\pi]$ and $\delta\in\left\{\Re(\delta)\leq 1, \left|\delta\right|\leq R\right\}$, where $R>1$, and noting that $\frac{\beta_{i}^{*}(\lambda(1-\delta))}{\theta+\lambda(1-\beta_{3}^{*}(\lambda(1-\delta)))}\leq \frac{1}{\theta}$,
\begin{displaymath}
\begin{array}{rl}
\left|k^{*}(\delta)\right|\leq&\left|\frac{\theta_{1}r_{1}}{\theta}+\frac{\theta_{2}r_{2}}{\theta}+\frac{2cos\phi\sqrt{\theta_{1}r_{1}\theta_{2}r_{2}}}{\theta}\right|\leq \frac{\theta_{1}r_{1}}{\theta}+\frac{\theta_{2}r_{2}}{\theta}+\frac{2\sqrt{\theta_{1}r_{1}\theta_{2}r_{2}}}{\theta}\\<&\frac{\theta_{1}r_{1}}{\theta}+\frac{\theta_{2}r_{2}}{\theta}+\frac{\theta_{1}r_{2}}{\theta}+\frac{\theta_{2}r_{1}}{\theta}<\frac{\theta_{1}+\theta_{2}}{\theta}=1<R=\left|\delta\right|.
\end{array}
\end{displaymath}
\paragraph{Proof of Lemma \ref{lem2}}We focus only in $L_{1}$. For $\delta=\delta(\phi)\in \mathbb{R}$,
\begin{equation}
\begin{array}{rl}
\delta(\phi)=&\frac{\theta_{1}r_{1}\beta_{1}^{*}(\lambda(1-\delta(\phi)))+\theta_{2}r_{2}\beta_{2}^{*}(\lambda(1-\delta(\phi)))}{\Delta(\phi)}\\&+\frac{2\cos\phi\sqrt{\theta_{1}r_{1}\theta_{2}r_{2}}\sqrt{\beta_{1}^{*}(\lambda(1-\delta(\phi)))\beta_{2}^{*}(\lambda(1-\delta(\phi)))}}{\Delta(\phi)},
\end{array}
\label{sdf}
\end{equation}
and $\delta(\phi) = \delta(\phi+\pi)$ for every $\phi\in[0,2\pi]$. Rewrite $w_{1}(\phi)$ as follows:
\begin{displaymath}
\begin{array}{rl}
w_{1}(\phi)=&[a(\delta(\phi))+b(\delta(\phi))\cos(\phi)]+ib(\delta(\phi))\sin(\phi),\\
a(\delta(\phi))=&\frac{2}{\Delta(\phi)}\theta_{1}r_{1}\beta_{1}^{*}(\lambda(1-\delta(\phi))),\\
b(\delta(\phi))=&\frac{2}{\Delta(\phi)}\sqrt{r_{1}\theta_{1}r_{2}\theta_{2}}\sqrt{\beta_{1}^{*}(\lambda(1-\delta(\phi)))\beta_{2}^{*}(\lambda(1-\delta(\phi)))}.
\end{array}
\end{displaymath}
Since $a(\delta(\phi))$ and $b(\delta(\phi))$ are the differentiable functions of $\delta$, we only need to show that $\delta(\phi)$ is a continuous differentiable function of $\phi$. By differentiating (\ref{sdf}) in $\phi$, we can show after some algebra that under the stability conditions $L_{1}$ is smooth and non-self intersecting. 
\paragraph{Proof of Lemma \ref{lem3}}
If $r_{2}\theta_{1}=r_{1}\theta_{2}$. Since $\delta(0)=1$, the contours
\begin{displaymath}
\begin{array}{c}
\frac{w_{j}(\phi)}{2r_{j}}-\widetilde{\beta}_{j}^{*}(w_{1}(\phi)/2r_{1},w_{2}(\phi)/2r_{2}),\,j=1,2,
\end{array}
\end{displaymath}
are smooth and have only two real points of which one is negative and the other that corresponds to $\delta(0)=1$ is located at zero, where the contours have vertical tangents. Therefore
\begin{displaymath}
\begin{array}{c}
ind_{\phi\in[0,2\pi]}[\frac{w_{2}(0)}{2r_{2}}-\widetilde{\beta}_{2}^{*}(w_{1}(0)/2r_{1},w_{2}(0)/2r_{2})]=-\frac{1}{2},\\
ind_{\phi\in[0,2\pi]}[\frac{w_{1}(0)}{2r_{1}}-\widetilde{\beta}_{1}^{*}(w_{1}(0)/2r_{1},w_{2}(0)/2r_{2})]=\frac{1}{2}.
\end{array}
\end{displaymath}

Thus, since $ind_{\phi\in[0,2\pi]}\frac{w_{1}(\phi)}{2r_{1}}=1$, $ind_{\phi\in[0,2\pi]}\frac{w_{2}(\phi)}{2r_{2}}=-1$, then $\chi=1-(-1)+(-1/2)-1/2=1$.\\
If $r_{2}\theta_{1}<r_{1}\theta_{2}$.
\begin{enumerate}
\item If $r_{1}\widehat{\theta}_{2}\geq r_{2}c^{2}/\widehat{\theta}_{1}$, then $w_{2}(0)/2r_{2}\geq 1$. In this case $ind_{\phi\in[0,2\pi]}[\frac{w_{2}(\phi)}{2r_{2}}-\widetilde{\beta}_{2}^{*}(w_{1}(0)/2r_{1},w_{2}(0)/2r_{2})]=-1$.
\item If $\theta_{2}r_{1}\widetilde{\beta}_{1}^{*}(w_{1}(0)/2r_{1},w_{2}(0)/2r_{2})>\Delta_{1}(0)r_{2}\beta_{2}^{*}(\lambda(1-\delta(0)))$, then we can guarantee that $\frac{w_{1}(0)}{2r_{1}}-\widetilde{\beta}_{1}^{*}(w_{1}(0)/2r_{1},w_{2}(0)/2r_{2})<0$. which in turn implies that $ind_{\phi\in[0,2\pi]}[\frac{w_{2}(\phi)}{2r_{2}}-\widetilde{\beta}_{2}^{*}(w_{1}(0)/2r_{1},w_{2}(0)/2r_{2})]=0$. In such a case, $\chi=1-(-1)+(-1)-0=1.$
\end{enumerate}

%
%

%


\end{document}